\documentclass{elsart}

\usepackage{amsmath, amsfonts}
\usepackage{amssymb}
\usepackage{stmaryrd}
\usepackage{graphicx, subfigure}
\usepackage{url}

\newcommand{\R}{{\mathbb R}}

\newcommand{\rmu}{{\rm i}}
\newcommand{\curl}{\mbox{\rm curl}\; }  

\newcommand{\Prond}{\mathcal{P}}

\newcommand{\Trond}{\mathcal{T}}

\newcommand{\jump}[1]{\llbracket #1 \rrbracket}
\newcommand{\average}[1]{\{#1\}}

\renewcommand{\vec}[1]{\mbox{\boldmath $#1$}} 

\renewcommand{\epsilon}{\varepsilon}
\renewcommand{\emph}{\textbf}

\renewcommand{\bar}{\overline}
\renewcommand{\conj}[1]{\overline{#1}}

\newtheorem{proposition}{{\sc Proposition}}

\newcommand{\QED}{\hspace*{\fill}\rule{2.5mm}{2.5mm}} 
\newenvironment{proof}{\noindent{\em Proof\ }}{\QED} 

\journal{Journal of Computational and Applied Mathematics}

\begin{document}
\begin{frontmatter}

  \title{Solution of the time-harmonic Maxwell equations using discontinuous
    Galerkin methods}

  \author[Unice,INRIA]{V. Dolean\corauthref{cor}},
  \corauth[cor]{Corresponding author.}
  \ead{dolean@math.unice.fr}
  \author[INRIA]{H. Fol},
  \ead{Hugo.Fol@inria.fr}
  \author[INRIA]{S. Lanteri},
  \ead{Stephane.Lanteri@inria.fr}
  \author[INRIA]{R. Perrussel}
  \ead{Ronan.Perrussel@inria.fr}

  \address[Unice]{Universit\'e de Nice Sophia-Antipolis, 
    Laboratoire J.-A. Dieudonn\'e, Parc Valrose, 
    06902 Nice Cedex 02}
  \address[INRIA]{INRIA Sophia-Antipolis, 2004 rte des Lucioles, 
    06902 Sophia-Antipolis Cedex}

  \begin{abstract}
    We present numerical results  concerning the solution of the time-harmonic
    Maxwell  equations  discretized  by  discontinuous Galerkin  methods.   In
    particular, a numerical study of the convergence, which compares different
    strategies proposed in the  literature for the elliptic Maxwell equations,
    is performed in the two-dimensional case.
  \end{abstract}

\end{frontmatter}

\section{Introduction}
\label{sec:intro}

This  work is  concerned  with  the numerical  solution  of the  time-harmonic
Maxwell   equations   discretized  by   discontinuous   Galerkin  methods   on
unstructured meshes.  Our motivation for using a discontinuous Galerkin method
is the  enhanced flexibility  compared to the  conforming edge  element method
\cite{monk03}:   for   instance,  dealing   with   non-conforming  meshes   is
straightforward  and  the choice  of  the  local  approximation space  is  not
constrained.  Nonetheless, before taking  full advantage of these features, it
is required to carefully study the basic ingredients of the method such as the
choice of  the numerical flux  at the interface between  neighboring elements.
In the  context of  time-harmonic problems, the  design of  efficient solution
strategies for  the resulting  sparse linear systems  is an  equally important
question.

Previous  works  have shown  convergence  results  for discontinuous  Galerkin
methods applied to the time-harmonic Maxwell equations, studied in the form of
second-order  vector  wave  equations.   Most  of  these  works  use  a  mixed
formulation \cite{perugia02,houston05b} but  discontinuous Galerkin methods on
the  non-mixed formulation  have recently  been proved  to  converge (interior
penalty technique \cite{houston05,buffa05} as  well as the local discontinuous
Galerkin method \cite{buffa05}).  The  convergence properties of these methods
in  the time-domain case  have been  studied in  \cite{fezoui05} when  using a
centered flux and  in \cite{hesthaven02} when using an  upwind flux.  The case
of  the   upwind  flux  has   been  analyzed  in   \cite{helluy-dayma:94}  and
\cite{helluy:94} for  the time-harmonic problems and the  convergence has been
proved only for  a perturbed problem.  The general  case of Friedrichs systems
and  the  elliptic  Maxwell  equations  in  particular  has  been  treated  in
\cite{ern06a}  and  \cite{ern06b}.   However,  to  our  knowledge,  no  direct
convergence      analysis      on      the      first-order      time-harmonic
system~\eqref{eq:maxwell_3d}  has  been  conducted  so far,  which  should  be
useful,   for  instance,   when   using  an   upwind   flux  (see   subsection
\ref{sec:flux_num}). The main  contribution of this work is  a numerical study
of the  convergence of  discontinuous Galerkin methods  based on  centered and
upwind fluxes applied  to the first-order time-harmonic Maxwell  system in the
two-dimensional case.

The paper  is organized as follows:  in section \ref{sec:form}  we present the
discretization method as  well as the different kind  of fluxes considered. In
section \ref{sec:convdis}, the convergence properties are recalled in the case
of  the  elliptic  Maxwell  equations  and the  solvability  of  the  discrete
perturbed problem  is analyzed in the  case of the centered  flux.  In section
\ref{sec:res_2d} the  numerical convergence is  studied and confronted  to the
theoretical convergence  order. Different  numerical fluxes are  then compared
on two distinct examples.

\section{Discretization of the first-order time-harmonic Maxwell system}
\label{sec:form}
\subsection{Formulation of the continuous problem}
\label{sec:cont_prob}

The  system of  non-dimensionalized time-harmonic  Maxwell's equations  can be
written in the following form:
\begin{equation}
  \label{eq:maxwell_3d}  
  \left\{
    \begin{aligned}
      \rmu \omega \epsilon_r \vec{E} - \curl \vec{H} & = -\vec{J}, \\
      \rmu \omega \mu_r \vec{H} + \curl \vec{E} & = 0,
    \end{aligned}
  \right.
\end{equation}
where $\vec{E}$ and $\vec{H}$ are the unknown electric and magnetic fields and
$\vec{J}$  is a  known  current source.   The  parameters $\varepsilon_r$  and
$\mu_r$ are  respectively the complex-valued  relative dielectric permittivity
(integrating   the   electric   conductivity)   and  the   relative   magnetic
permeability;  we consider  here  the  case of  linear  isotropic media.   The
angular  frequency   of  the   problem  is  given   by  $\omega$.    We  solve
Equa\-tions~\eqref{eq:maxwell_3d}  in a  bounded domain  $\Omega$, and  on its
boundary  $\partial  \Omega  =  \Gamma_{a}  \cup \Gamma_{m}$,  we  impose  the
following boundary conditions:
\begin{equation} 
  \label{eq:condlim_3d}
  \begin{aligned}
    & \text{- a perfect electric conductor condition on } \Gamma_m,
    \text{ {\sl ie}: } \vec{n} \times \vec{E} = 0 \text{ on }
    \Gamma_m, \\
    & \text{- a Silver-M\"uller (first-order absorbing boundary)
      condition on } \Gamma^a, \text{ {\sl ie}: } \\
    & \vec{n} \times \vec{E} + \vec{n} \times (\vec{n} \times \vec{H} 
    ) = \vec{n} \times \vec{E}^{\rm inc} + \vec{n} \times
    (\vec{n} \times \vec{H}^{\rm inc} ) \text{ on } \Gamma_a.
  \end{aligned}
\end{equation}
The  vectors  $\vec{E}^{\rm  inc}$   and  $\vec{H}^{\rm  inc}$  represent  the
components  of  an  incident  electromagnetic  wave. We  can  further  rewrite
\eqref{eq:maxwell_3d}+\eqref{eq:condlim_3d},  assuming $\vec  J$ equals  to 0,
under the following form:
\begin{equation}
  \label{eq:3Dformfort}
  \left\{
    \begin{aligned}
      & \rmu \omega G_0 \vec{W} + G_x \partial_x \vec{W} + G_y
      \partial_y \vec{W} + G_z
      \partial_z \vec{W} = 0 \text{ in } \Omega, \\
      & (M_{\Gamma_m} - G_{\vec{n}}) \vec{W} = 0 \text{ on } \Gamma_m, \\
      & (M_{\Gamma_a} - G_{\vec{n}}) (\vec{W} - \vec{W}_{\rm inc}) = 0
      \text{ on } \Gamma_a.
    \end{aligned}
  \right.
\end{equation}
where $\vec{W} = \begin{pmatrix} \vec{E}  \\ \vec{H} \end{pmatrix}$ is the new
unknown vector and $ G_0 =
\begin{pmatrix}
  \epsilon_r I_3 & 0_{3 \times 3} \\
  0_{3 \times 3} & \mu_r I_3
\end{pmatrix}.  $ Denoting by $(\vec e^x,  \vec e^y, \vec e^z )$ the canonical
basis of $\R^3$, the matrices $G_l$ with $l\in\{x, y, z\}$ are given by:
\begin{equation*}
  G_{l} =
  \begin{pmatrix}
    0_{3 \times 3} & N_{\vec{e}^l} \\ N_{\vec{e}^l}^t & 0_{3 \times 3}
  \end{pmatrix}
  \ \text{where for a vector } \vec n, \ N_{\vec{n}} =
  \begin{pmatrix}
    0 & \vec{n}_z & -\vec{n}_y \\ -\vec{n}_z & 0 & \vec{n}_x \\ \vec{n}_y &
    -\vec{n}_x & 0
  \end{pmatrix}.
\end{equation*}
In the  following we denote  by $G_{\vec{n}}$ the  sum $G_x \vec  n_x+G_y \vec
n_y+G_z \vec n_z$ and by  $G_{\vec{n}}^+$ and $G_{\vec{n}}^-$ its positive and
negative   parts\footnote{If  $G_{\vec  n}   =  T   \Lambda  T^{-1}$   is  the
  eigenfactorization  then  $G^\pm_{\vec  n}  = T  \Lambda^\pm  T^{-1}$  where
  $\Lambda^+$ (resp.  $\Lambda^-$) only  gathers the positive (resp. negative)
  eigenvalues.}.  We also define  $\lvert G_{\vec{n}} \lvert = G_{\vec{n}}^+ -
G_{\vec{n}}^-$.  In  order to take  into account the boundary  conditions, the
matrices $M_{\Gamma_m}$ and $M_{\Gamma_a}$ are given by:
\begin{equation*}
  M_{\Gamma_m} =
  \begin{pmatrix}
    0_{3 \times 3} & N_{\vec{n}} \\ -N_{\vec{n}}^t & 0_{3 \times 3}
  \end{pmatrix}
  \text{ and } M_{\Gamma_a} = \lvert G_{\vec{n}} \lvert.
\end{equation*}
See \cite{dolean_05} for further details on the derivation of this formulation.

\subsection{Discretization}
\label{sec:discret}

Let $\Omega_h$ denote a discretization of  the domain $\Omega$ into a union of
conforming elements (tetrahedral or hexahedral elements)
$$\displaystyle\overline{\Omega}_h  =
\bigcup_{K  \in  \Trond_h}  K.
$$
We look for the approximate solutions $\vec W_h = \begin{pmatrix} \vec{E}_h \\
  \vec{H}_h \end{pmatrix}$ of \eqref{eq:3Dformfort}  in $V_h \times V_h$ where
the function space $V_h$ is defined by:
\begin{equation}
  \label{eq:def_Vh} 
  V_h  = \left\{  \vec{V}  \in  [L^2(\Omega)]^3 \  /  \  \forall K  \in
    \Trond_h, \ \ \vec{V}_{\lvert K} \in \Prond(K) \right\}.
\end{equation}
The term  $\Prond(K)$ denotes a space  of polynomial functions  on the element
$K$.    We   take   the   scalar    product   of   the   first   equation   of
(\ref{eq:3Dformfort}) by  a sufficiently smooth vector field  $\vec{V}$ and we
integrate over an element $K$ of the mesh $\Trond_h$:
$$
\int_K \rmu \omega \left( G_0 \vec{W} \right)^t \bar{\vec{V}} dx + \int_K
\left( \sum_{l \in \{x,y,z\}} G_l \partial_l \vec{W} \right)^t \bar{\vec{V}}
dx = 0.
$$ 
By using  Green's formula  we obtain a  weak formulation involving  a boundary
term.  This term  is replaced in discontinuous Galerkin  methods by a function
$\Phi_{\partial K}$ which is usually  referred as the numerical flux (see also
Ern and  Guermond \cite{ern06a,ern06b});  the aim is  then to  determine $\vec
W_h$ in $V_h \times V_h$ such that:
\begin{equation}
  \label{eq:weak_form} 
  \begin{aligned}
    \int_K \rmu \omega \left( G_0 \vec{W}_h \right)^t \bar{\vec{V}} dx -
    \int_K \vec{W}_h^t \left( \sum_{l \in \{x, y, z\}} G_l \partial_l
      \bar{\vec{V}} \right) dx + \int_{\partial K} \left( \Phi_{\partial
        K}(\vec{W}_h) \right)^t \bar{\vec{V}} = 0,&
    \\
    \forall \vec V \in V_h \times V_h.&
  \end{aligned}
\end{equation}
In  order to  couple  the element  $K$  with its  neighbors  for ensuring  the
consistency of the  discretization, this numerical flux can  be defined in the
following way:
\begin{equation}
  \Phi_{\partial K}(\vec W_h) = \left\{
    \begin{aligned}
      & I_{FK} S_{F} \jump{\vec{W}_h} + I_{FK} G_{\vec{n}_F} \average{\vec{W}_h}
      \text{ if } F \in \Gamma^0,
      \\
      & \frac{1}{2}(M_{F, K} + I_{FK} G_{\vec{n}_F})\vec{W}_h \text{ if } F \in
      \Gamma^m,
      \\
      & \frac{1}{2}(M_{F, K} + I_{FK} G_{\vec{n}_F})\vec{W}_h -
      \frac{1}{2}(M_{F, K} - I_{FK} G_{\vec{n}_F}) \vec{W}^{\rm inc} \text{ if
      } F \in \Gamma^a,
    \end{aligned} \right.
\end{equation}
where  $\Gamma^0$, $\Gamma^a$ and  $\Gamma^m$ respectively  denote the  set of
interior  faces, the  set  of faces  on $\Gamma_a$  and  the set  of faces  on
$\Gamma_m$.  $I_{FK}$  stands for the incidence matrix  between oriented faces
and elements whose entries are given by:
\begin{equation*} I_{FK} = \left\{
    \begin{aligned}
      0 & \text{ if the face $F$ does not belong to element $K$,}
      \\
      1 & \text{ if $F \in K$ and their orientations match,}
      \\
      -1 & \text{ if $F \in K$ and their orientations do not match.}
    \end{aligned} \right.
\end{equation*}
We also define respectively the jump and the average of $\vec V$ on a face $F$
shared by two elements $K$ and $\tilde K$:
\begin{equation*}
  \jump{\vec{V}}  = I_{FK} \vec{V}_K  + I_{F \tilde K}  \vec{V}_{\tilde K}
  \text{ and } \average{\vec{V}} = \frac{1}{2} (\vec{V}_K + \vec{V}_{\tilde K}).
\end{equation*}
Finally, the matrix  $S_F$ allows to penalize  the jump of a field  or of some
components of this given field on the  face $F$ and the matrix $M_{F,K}$ to be
defined later insures the  asymptotic consistency with the boundary conditions
of the continuous problem.

\subsection{Choice of the numerical flux}
\label{sec:flux_num}
In this study, we aim at comparing the properties of three classical numerical
fluxes:

\noindent -~ {\bf  a centered flux} (see \cite{fezoui05}  for the time-domain equivalent).
In this case $S_{F} = 0$ for all the faces $F$ and, for the boundary faces, we use:
\begin{equation*} M_{F, K} = \left\{
    \begin{aligned} & I_{FK}
      \begin{pmatrix} 0_{3 \times 3} & N_{\vec{n}_F} \\ -N_{\vec{n}_F}^t &
        0_{3 \times 3}
      \end{pmatrix} \text{ if } F \in \Gamma^m, \\ & \lvert G_{\vec{n}_F}
      \lvert \text{ if } F \in \Gamma^a.
    \end{aligned} \right.
\end{equation*}

\noindent -~ {\bf an upwind flux} (see \cite{ern06a,piperno00}). In this case:
\begin{equation*} S_F =
  \begin{pmatrix}
    \alpha_F^E N_{\vec{n}} N_{\vec{n}}^t & 0_{3 \times 3} \\ 0_{3 \times 3} &
    \alpha_F^H N_{\vec{n}}^t N_{\vec{n}}
\end{pmatrix}, \ M_{F, K} =
\begin{pmatrix}
  \eta_F N_{\vec{n}_F} N_{\vec{n}_F}^t & I_{FK} N_{\vec{n}_F} \\ -I_{FK}
  N_{\vec{n}_F}^t & 0_{3 \times 3}
  \end{pmatrix} \ \forall F \in \Gamma^m,
\end{equation*}
with $\alpha_F^E$,  $\alpha_F^H$ and $\eta_F$ equals to  $1/2$ for homogeneous
media.  The definition  of $M_{FK}$ for $F$ in $\Gamma^a$  is identical to the
centered case.

\noindent  -~ {\bf  a partially  penalized upwind  flux}  (local Discontinuous
Galerkin  method, see  \cite{cockburn98}).  This  flux is  characterized  by a
penalization coefficient given by:
\begin{equation*} S_F = \tau_F h_F^{-1}
  \begin{pmatrix} N_{\vec{n}_F} N_{\vec{n}_F}^t & 0 \\ 0 & 0
  \end{pmatrix}, \ M_{F, K} =
  \begin{pmatrix}
    \eta_F h_F^{-1} N_{\vec{n}_F}  N_{\vec{n}_F}^t & I_{FK} N_{\vec{n}_F} \\
    -I_{FK} N_{\vec{n}_F}^t & 0_{3 \times 3}
  \end{pmatrix} \ \forall F \in \Gamma^m.
\end{equation*}
The definition  of $M_{FK}$  for $F$  in $\Gamma^a$ is  also identical  to the
centered case.

\section{Convergence properties of the discretized problem}
\label{sec:convdis}

We are interested  in assessing these numerical fluxes  for the discretization
of \eqref{eq:3Dformfort}.   Firstly, we  want the best  asymptotic convergence
order in $L^2$-norm for the electric and magnetic field for a fixed polynomial
order  approximation on an  unstructured mesh.  Secondly, a  minimal numerical
dispersion  is also  needed.  In  the  following we  will focus  on the  first
criterion.  The  asymptotic convergence order in $L^2$-norm  between the exact
solution  $(\vec{E},  \vec{H})$  and  the  approximate  solution  $(\vec{E}_h,
\vec{H}_h)$ corresponds to the  largest real coefficients $\beta$ and $\gamma$
such that:
\begin{equation}
  \exists C_1, C_2,  h_0 > 0, \ \forall h > h_0,  \ \| \vec{E} - \vec{E}_h
  \|_{L^2(\Omega)} \leq  C_1 h^\beta  \text{ and }  \| \vec{H} -  \vec{H}_h \|_{L^2(\Omega)}
  \leq C_2 h^\gamma,
  \label{eq:gamma_beta}
\end{equation}
where  $h$ is  the mesh  size.  Let  us note  that in  the  numerical examples
proposed in Section \ref{sec:res_2d},  we will often equivalently consider the
evolution of the  norm of the error  against the square root of  the number of
degrees  of  freedom (dofs),  in  order  to  deduce coefficients  $\beta$  and
$\gamma$.

We first  recall in Table~\ref{tab:Max_ell} below  the theoretical convergence
order  for   the  elliptic  Maxwell  equations   \cite{ern06a,ern06b},  for  a
sufficiently smooth solution and when  the local function space $\Prond(K)$ is
$[P_k(K)]^3$  \textsl{i.e.}   the  space   of  vectors  whose  components  are
polynomials of order at most $k$.   When using the flux with a penalization of
$\vec E$, similar convergence results are proved for the time-harmonic Maxwell
equations in \cite{buffa05}.
\begin{table}[htbp] \centering
  \begin{tabular}[c]{c|c|c|c}  flux &  centered &  upwind &  penalization of  $\vec{E}$ \\
\hline field $\vec{E}$ & $k$ & $k+1/2$ & $k+1$ \\ \hline field $\vec{H}$ & $k$ & $k+1/2$ &
$k$
  \end{tabular}
  \caption{Theoretical convergence order for the elliptic Maxwell equations.}
  \label{tab:Max_ell}
\end{table}

\subsection{Solution of the discretized perturbed problem}
\label{sec:resdis}

A few comments need to be stated concerning the convergence properties of such
a scheme applied  to the first-order formulation of  the time-harmonic Maxwell
equations.  First  of all, the  case of the  upwind flux has been  analyzed in
\cite{helluy-dayma:94}  for  the  perturbed  Maxwell  problem,  that  is  when
$i\omega$  is  replaced  by  $\nu+i\omega$  with  $\nu$  a  strictly  positive
parameter. For a  sufficiently regular solution the norm  of the error behaves
as $h^{p+1/2}$ where $h$ is the mesh parameter.

The case  of the  centered flux  has been studied  in \cite{fezoui05}  for the
time-domain Maxwell equations  and in this case the norm  of the error behaves
as $h^{p}$ where $h$ is the mesh parameter. For the time-harmonic equations no
convergence  proofs  are  available  so  far.   We can  only  study  here  the
solvability of  the discrete problem  in the case  of a perturbed  problem (we
replace $\rmu \omega$ by $\rmu \omega + \nu$ with $\nu > 0$) following an idea
used by Helluy \cite{helluy:94} in the case of the upwind flux. In the case of
the  perturbed  problem  and  assuming homogeneous  boundary  conditions,  the
formulation can be simply written as:
\begin{equation}
  \label{eq:formdirecthom}
  \left\{
  \begin{aligned}
    & \text{Find $\vec{W}_h$ in $V_h\times V_h$ such that:} \\
    & a(\vec{W}_h, \vec{V}) + b(\vec{W}_h, \vec{V}) = 0 , \ 
      \forall \vec{V} \in V_h \times V_h,
  \end{aligned}
  \right .
\end{equation}
with, $\forall \vec{U}, \vec{V} \in V_h \times V_h$: 
\begin{equation}
  \label{eq:forma}
  \begin{aligned}
    a(\vec{U}, \vec{V}) & = 
      \int_{\Omega_h} \left ((\rmu\omega + \nu)G_0\vec{U}\right )^t
      \conj{\vec{V}} dv + 
      \sum_{F \in \Gamma^a} 
      \int_F \left (\frac{1}{2}\lvert G_{\vec{n}_F} \lvert \vec{U}\right )^t 
      \conj{\vec{V}} ds \\
    & + \sum_{F \in \Gamma^m} 
      \int_F \left (\frac{1}{2}M_{F, K}\vec{U}\right )^t 
      \conj{\vec{V}} ds + 
      \sum_{F \in \Gamma^0} 
      \int_F \left (S_F\jump{\vec{U}}\right )^t 
      \jump{\conj{\vec{V}}}_F ds,
  \end{aligned}
\end{equation}
and:
\begin{equation}
  \label{eq:formb}
  \begin{aligned}
    b(\vec{U}, \vec{V}) & = 
      \sum_{K \in \Trond_h} 
      \int_K \left (\sum_{l \in \{x, y, z\}} G_l\partial_l (\vec{U})\right )^t 
      \conj{\vec{V}} dv \\ 
    & - \sum_{F \in \Gamma^a \cup \Gamma^m} 
      \int_F \left (\frac{1}{2}I_{FK}G_{\vec{n}_F}\vec{U}\right )^t 
      \conj{\vec{V}} ds \\
    & - \sum_{F \in \Gamma^0} 
      \int_F \left (G_{\vec{n}_F}\jump{\vec{U}}\right )^t 
      \average{\conj{\vec{V}}} ds.
  \end{aligned}
\end{equation}
We have the following result:
\begin{proposition}
The solution of problem \eqref{eq:formdirecthom} is null. 
\end{proposition}
\begin{proof}
First, considering the  fact that the  matrices $\lvert G_{\vec{n}_F} \lvert$,
$S_F$, $\Re  (G_0)$ and $-\Im (G_0)$  are hermitian and denoting  by $\mathcal
H(M_{F,K})$ the hermitian part of $M_{F,K}$ for $F$  in $\Gamma^m$, which is
equal to $\begin{pmatrix} 
\eta_F  N_{\vec{n}_F} N_{\vec{n}_F}^t & & 0_{3\times 3} \\ 
0_{3\times 3} & & 0_{3 \times 3} 
\end{pmatrix}$, one has:
\begin{equation}
  \label{eq:partiereel}
  \begin{aligned}
    \Re ( a(\vec{W}_h, \vec{W}_h) ) & =
    \int_{\Omega_h} \left ((\nu\Re (G_0) - \omega\Im (G_0))\vec{W}_h \right )^t
    \conj{\vec{W}_h} dv \\ & + 
    \sum_{F \in \Gamma^0} 
    \int_F \left (S_F\jump{\vec{W}_h}\right )^t 
    \jump{\conj{\vec{W}_h}}_F ds \\ & + 
    \sum_{F \in \Gamma^a} 
    \int_F \left (\frac{1}{2}\lvert G_{\vec{n}_F} \lvert\vec{W}_h \right )^t 
    \conj{\vec{W}_h} ds \\ & + 
    \sum_{F \in \Gamma^m} 
    \int_F \left (\frac{1}{2}\mathcal H (M_{F,K})\vec{W}_h\right )^t 
    \conj{\vec{W}_h} ds.
  \end{aligned}
\end{equation}
Then, we rewrite  using   the  corresponding  Green identity  an equivalent
expression of the sesquilinear form $b$:
\begin{equation}
  \label{eq:partiereel2}
  \begin{aligned}
    b(\vec{U}, \vec{V}) & = 
    - \sum_{K \in \Trond_h}
    \left [\int_K\vec{U}^t
           \left (\sum_{l \in \{x, y, z\}} G_l\partial_l (\conj{\vec{V}})\right ) dv -
    \right . \\ & \left . \hspace{10ex} 
    \sum_{F \in \partial K} 
    \int_F (I_{FK}G_{\vec{n}_F}\vec{U}_{\lvert K})^t \vec{V}_{\lvert K} ds 
    \right ] \\ & - 
    \sum_{F \in \Gamma^a} 
    \int_F \left (\frac{1}{2} I_{FK}G_{\vec{n}_F}\vec{U}\right )^t
    \conj{\vec{V}} ds \\ & -
    \sum_{F \in \Gamma^0} 
    \int_F \left (G_{\vec{n}_F}\jump{\vec{U}}\right )^t
    \average{\conj{\vec{V}}} ds, 
    \ \forall \vec{U}, \vec{V} \in V_h \times V_h.
  \end{aligned}
\end{equation}
By  noticing that on a  face $F \in \Gamma^0$  separating two elements $K$ and
$\tilde K$:
\begin{equation*}
  \begin{aligned}
  (G_{\vec{n}_F}\average{\vec{U}})^t\jump{\vec{V}} + 
  (G_{\vec{n}_F}\jump{\vec{U}})^t\average{\vec{V}} & =
  (I_{FK} G_{\vec{n}_F}\vec{U}_{\lvert K})^t\vec{V}_{\lvert K} \\ & + 
  (I_{F\tilde K} G_{\vec{n}_F}\vec{U}_{\lvert\tilde K})^t\vec{V}_{\lvert\tilde K},
  \end{aligned}
\end{equation*}
which is  in part  due to   the fact  that  $G_{\vec{n}_F}$ is  hermitian, one
deduces:
\begin{equation}
  \begin{aligned}
  b(\vec{U}, \vec{V}) & = 
  - \sum_{K \in \Trond_h} 
  \int_K\vec{U}^t
  \left (\sum_{l \in \{x, y, z\}}G_l\partial_l (\conj{\vec{V}})\right ) dv \\ & + 
  \sum_{F \in \Gamma^a} 
  \int_F \left (\frac{1}{2}I_{FK} G_{\vec{n}_F}\vec{U}\right )^t\conj{\vec{V}} ds \\ & + 
  \sum_{F \in \Gamma^0} 
  \int_F \left (G_{\vec{n}_F}\average{\vec{U}}\right )^t
  \jump{\conj{\vec{V}}}ds, 
  \ \forall \vec U, \vec{V} \in V_h \times V_h.
  \end{aligned}
\end{equation}  
Thus, it is now straightforward to see that $b$ is anti-hermitian and consequently:
\begin{equation*}
  \begin{aligned}
  \Re ( a(\vec{W}_h, \vec{W}_h) + b(\vec{W}_h, \vec{W}_h) ) & =
  \int_{\Omega_h}\left ((\nu\Re (G_0) - \omega\Im (G_0))\vec{W}_h\right )^t
  \conj{\vec{W}_h} dv \\ & + 
  \sum_{F \in \Gamma^0} 
  \int_F\left (S_F\jump{\vec{W}_h}\right )^t
  \jump{\conj{\vec{W}_h}}_F ds \\ & + 
  \sum_{F \in \Gamma^a} 
  \int_F\left (\frac{1}{2}\lvert G_{\vec{n}_F}\lvert\vec{W}_h\right )^t 
  \conj{\vec{W}_h} ds \\ & + 
  \sum_{F \in \Gamma^m} 
  \int_F \left (\frac{1}{2}\mathcal H (M_{F,K})\vec{W}_h\right )^t
  \conj{\vec{W}_h} ds,
  \end{aligned}
\end{equation*}
From \eqref{eq:formdirecthom}, $\Re ( a(\vec W_h, \vec W_h) + b(\vec W_h, \vec
W_h) )$  is also equal  to zero.  As  $\nu  \Re (G_0)  -  \omega \Im (G_0)$ is
positive definite  and  $\lvert G_{\vec  n_F}  \lvert$, $S_F$ and $\mathcal  H
(M_{F,K})$ are positive, the vector field $\vec{W}_h$ is null.
\end{proof}

\section{Numerical results}
\label{sec:res_2d}
In the  first part of this section  we will present a  numerical comparison of
different fluxes for a very simple  test case and different kind of meshes. In
the second part,  the results on a less trivial problem  are compared to those
obtained with the plane wave example.

\begin{figure}[!hbtp] \centering \subfigure[Initial mesh.]{
    \includegraphics[width=0.3\linewidth]{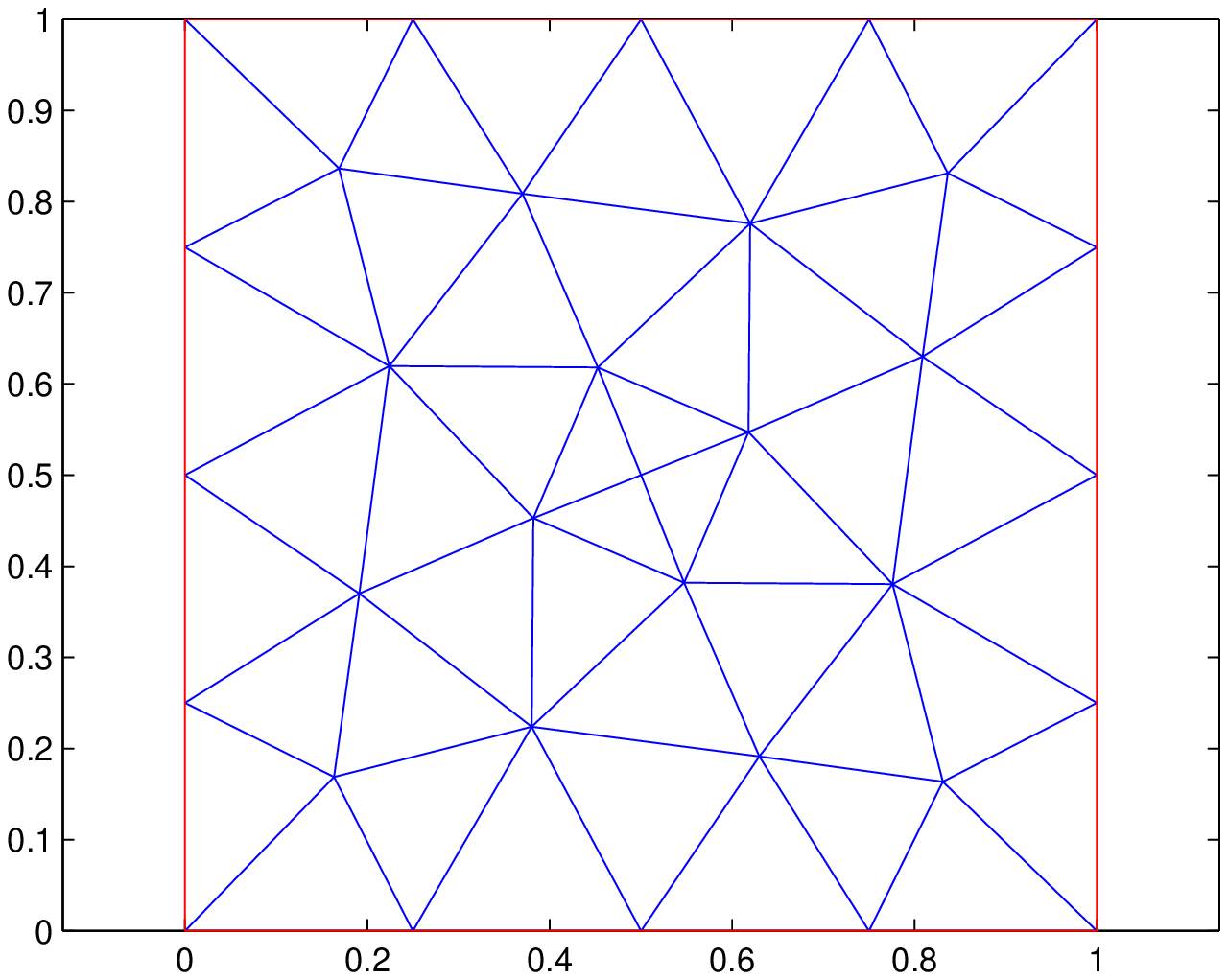}
    \label{fig:maill_init}} \hfill \subfigure[First refinement.]{
    \includegraphics[width=0.3\linewidth]{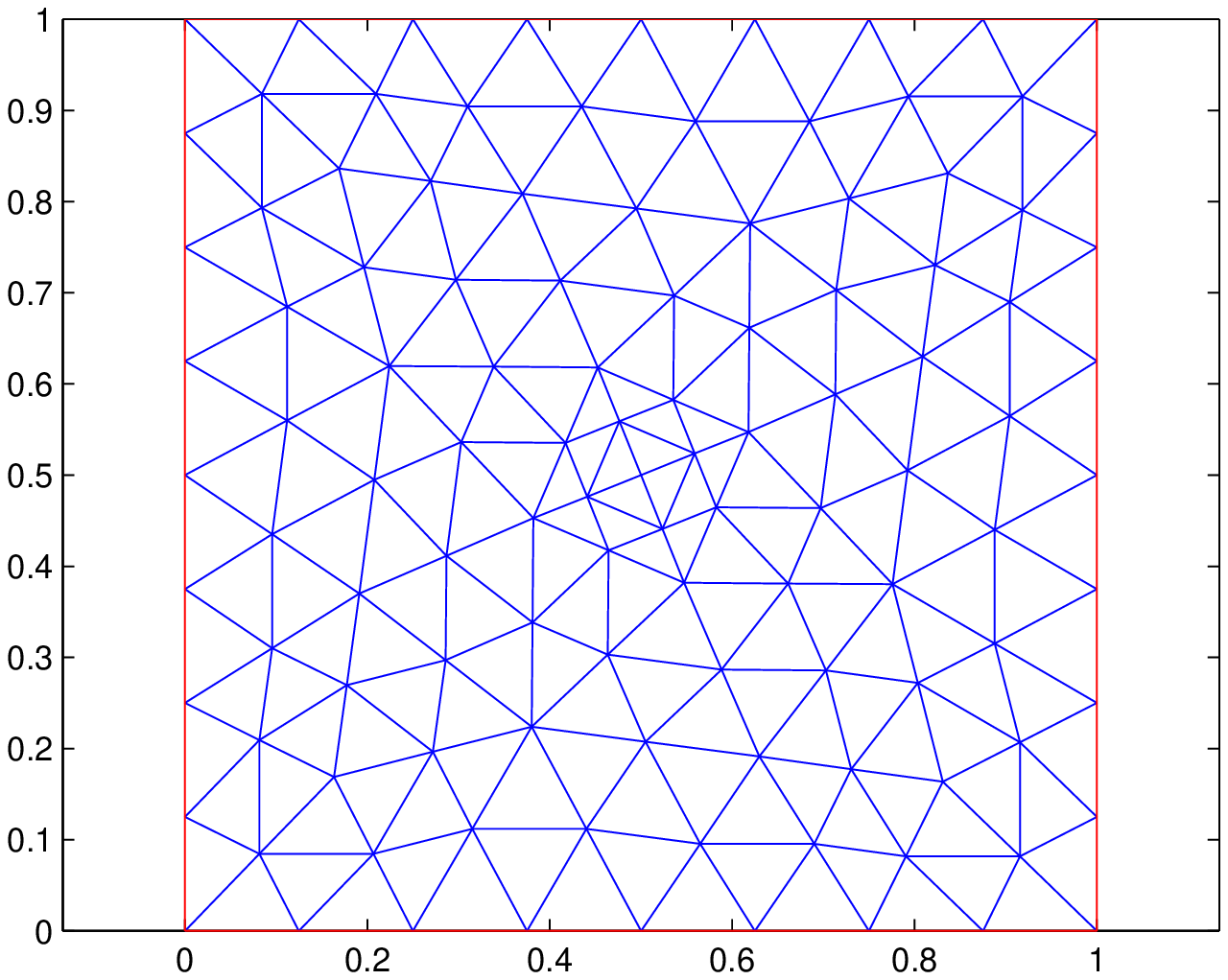}
    \label{fig:maill_1raff}} \hfill \subfigure[Second refinement.]{
    \includegraphics[width=0.3\linewidth]{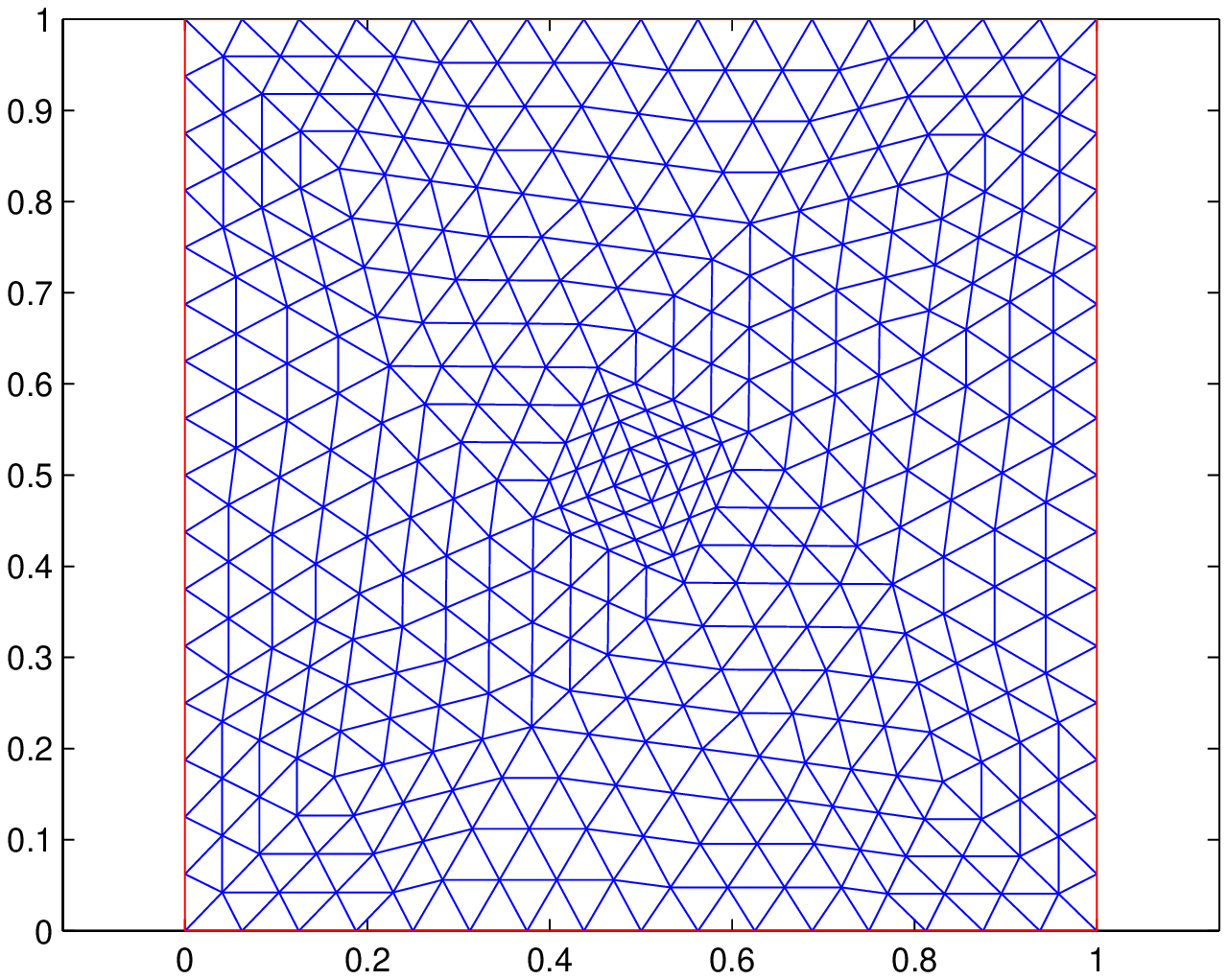}
    \label{fig:maill_2raff}}
  \caption{Initial mesh of the unit square and two uniform refinements.}
  \label{fig:maill_ustruct}
\end{figure}

We consider the case of an electric  transverse wave in the plane $(O, x, y)$.
In this case the components $\vec{E}_z$, $\vec{H}_x$ and $\vec{H}_y$ are zero.
We numerically  simulate the propagation of  a plane wave in  vacuum where the
incident    wave    is    given    by    $(\vec{E}^{\rm    inc}_x,\vec{E}^{\rm
  inc}_y,\vec{H}^{\rm inc}_z)=\exp (-  i \omega x)(0,1,1)$.  The computational
domain is the  unit square $\Omega = ]0; 1[^2$  and a Silver-M\"uller boundary
condition is imposed on the whole boundary, that is $\Gamma_a =
\partial \Omega$  and $\Gamma_m =\emptyset$.  The  parameters $\epsilon_r$ and
$\mu_r$ are set to $1$ everywhere and we choose $\omega=2\pi$.  We numerically
estimate the  asymptotic convergence  order of discontinuous  Galerkin methods
for the above problem using two different sequences of triangular meshes:

\noindent   -~  {\bf   uniformly   refined  meshes.}    The   first  mesh   of
Figure~\ref{fig:maill_init} is  uniformly refined  resulting in the  meshes of
Figures~\ref{fig:maill_1raff} and \ref{fig:maill_2raff}.

\begin{figure}[!hbtp] \centering \subfigure[$h = 1/8$.]{
    \includegraphics[width=0.3\linewidth]{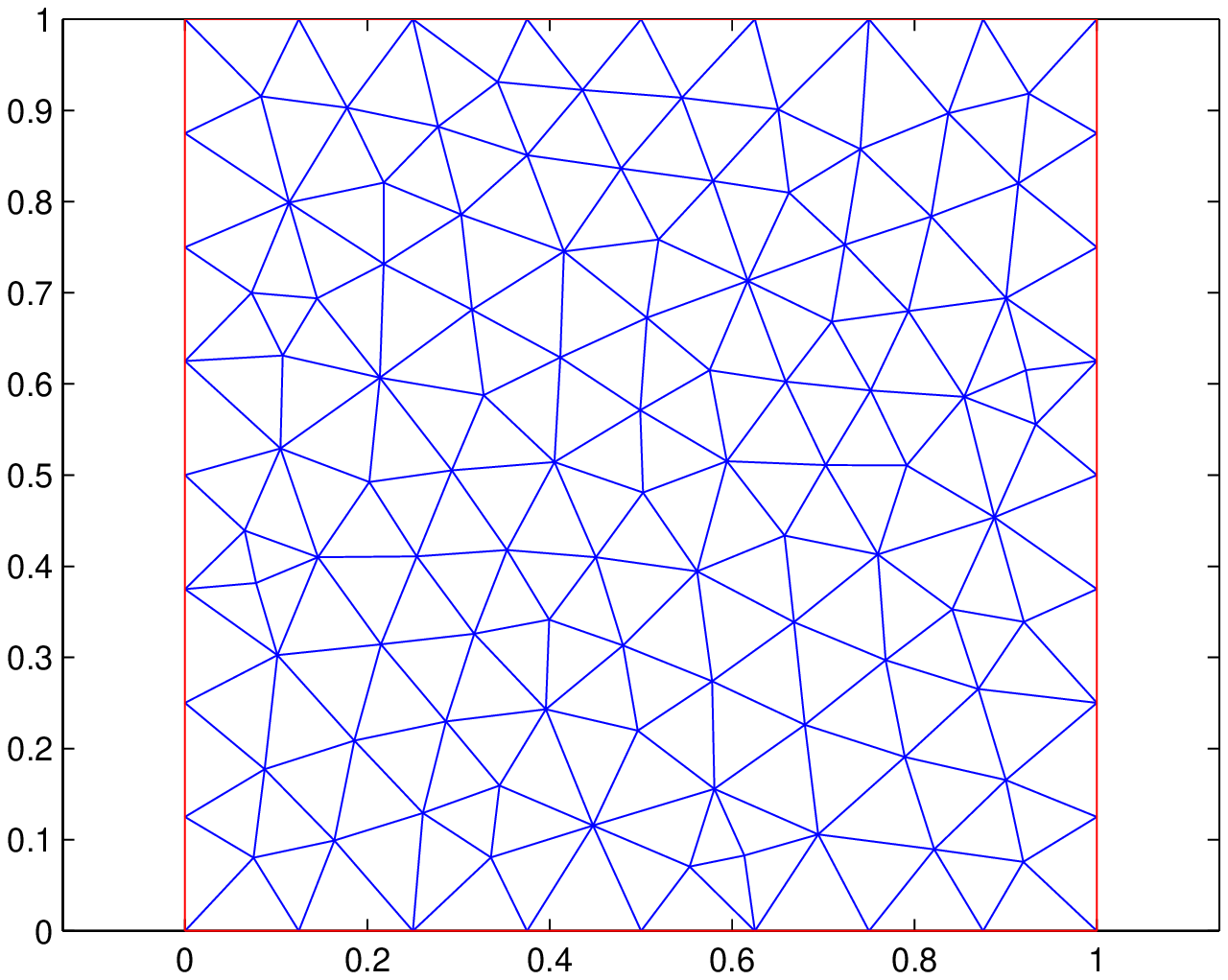}
    \label{fig:maill_0125}} \hfill \subfigure[$h= 1/16$.]{
    \includegraphics[width=0.3\linewidth]{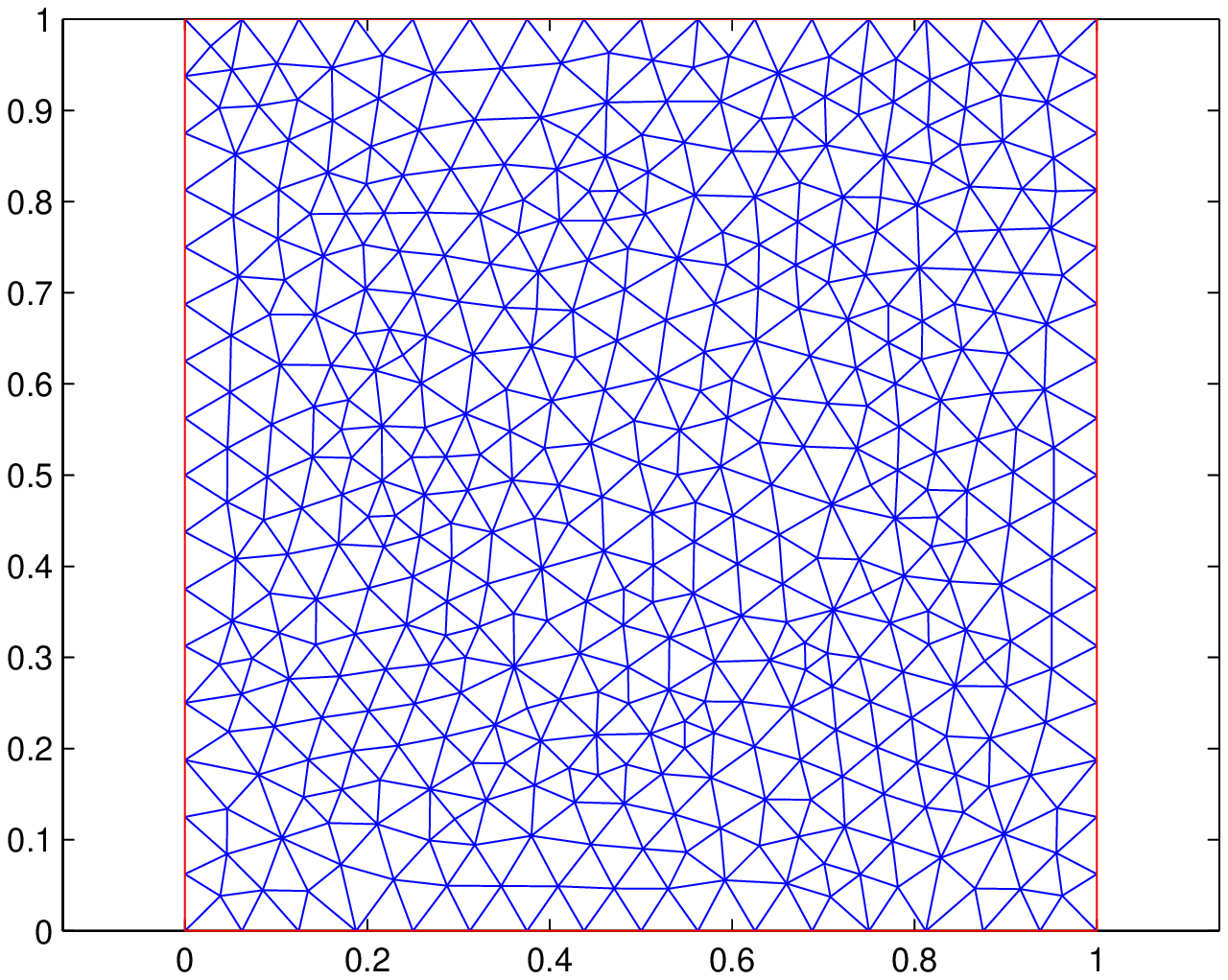}
    \label{fig:maill_00625}} \hfill \subfigure[$h= 1/32$.]{
    \includegraphics[width=0.3\linewidth]{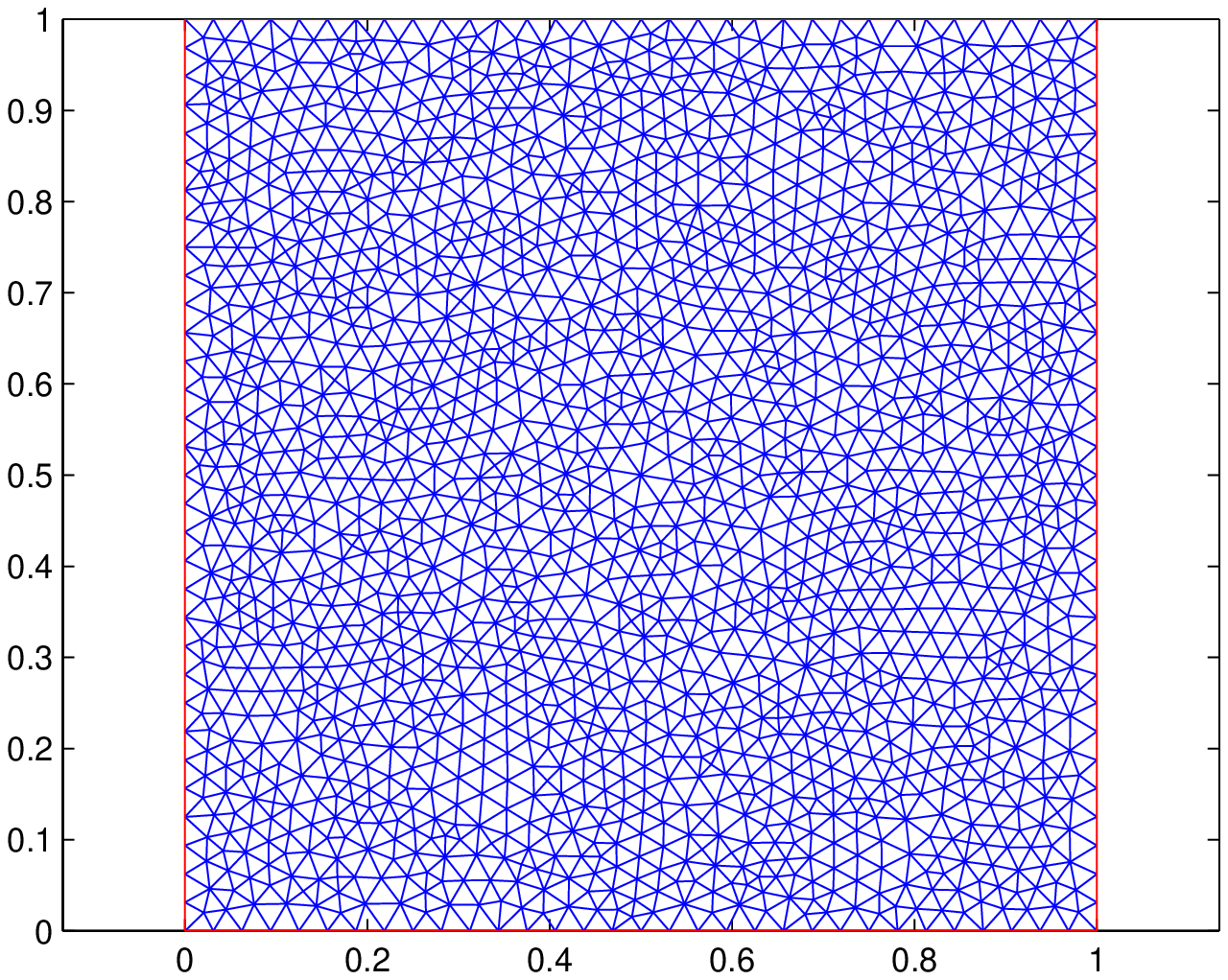}
    \label{fig:maill_003125}}
  \caption{First three independent unstructured meshes.}
  \label{fig:maill_tustruct}
\end{figure}

\noindent   -~   {\bf  independent   meshes.}    We   use  four   unstructured
(quasi-uni\-form)  independent meshes with  an imposed  maximal mesh  size $h$
(see Figure~\ref{fig:maill_tustruct} for the first three meshes). These meshes
are denoted by  $\Trond_i$ for $i=1,\ldots,4$ with $h$  in a decreasing order.
Thus $\Trond_{i+1}$ is not a refinement of $\Trond_i$.

Our implementation of  high order discontinuous Galerkin methods  makes use of
nodal basis functions with equi-spaced nodes.

\subsection{Convergence behavior using meshes obtained by uniform refinement}
\label{sec:ustruct}

\begin{figure}[!hbtp]
  \centering%
  \subfigure[$\| \vec{H} - \vec{H}_h \|_{L^2}$ against the square root of the
  number of dofs.]  {%
    \includegraphics[width=0.47\linewidth]{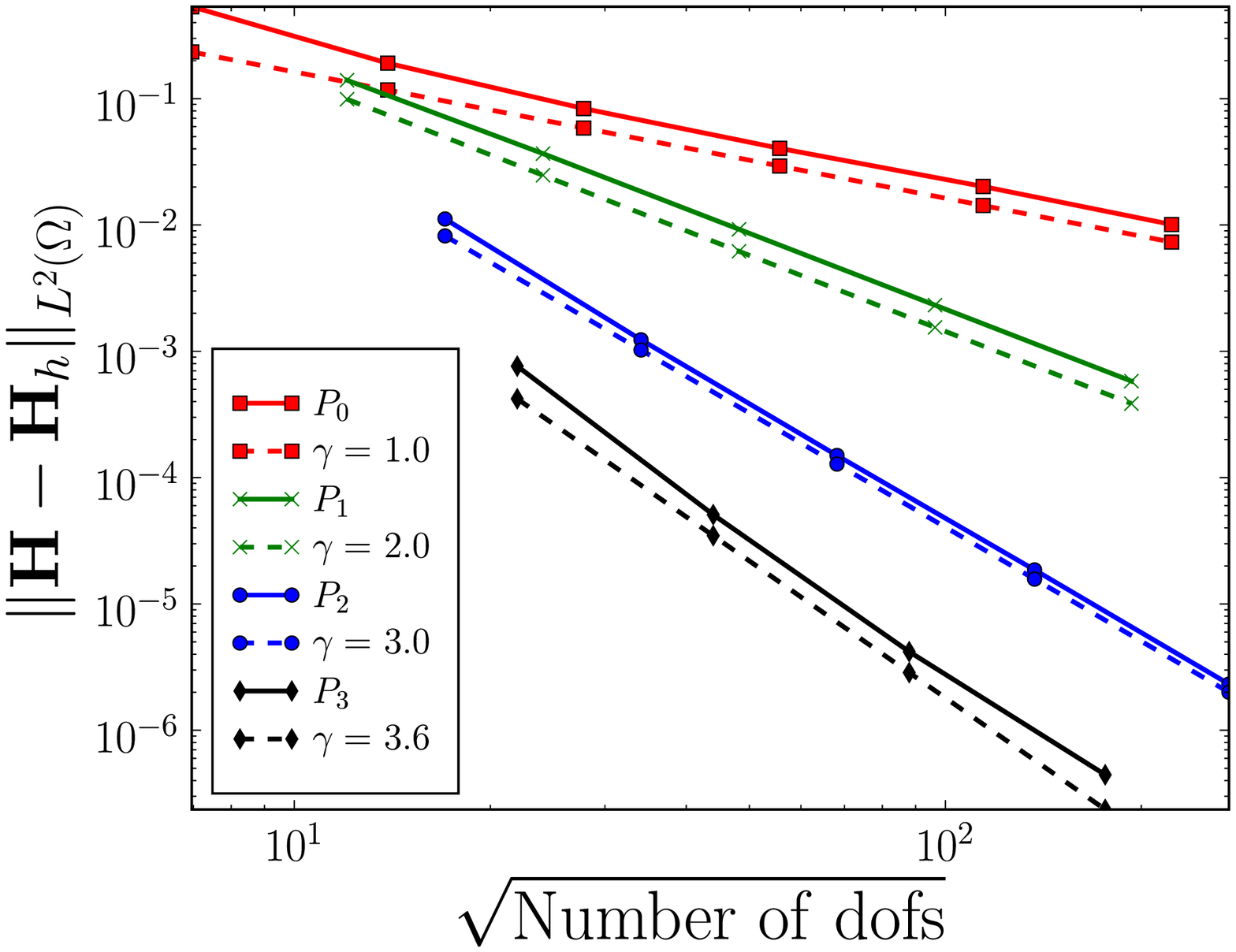} %
    \label{fig:res_c_H_pas} } %
  \hfill \subfigure[$\| \vec{E} - \vec{E}_h \|_{L^2}$ against the square root
  of the number of dofs.]  {%
    \includegraphics[width=0.47\linewidth]{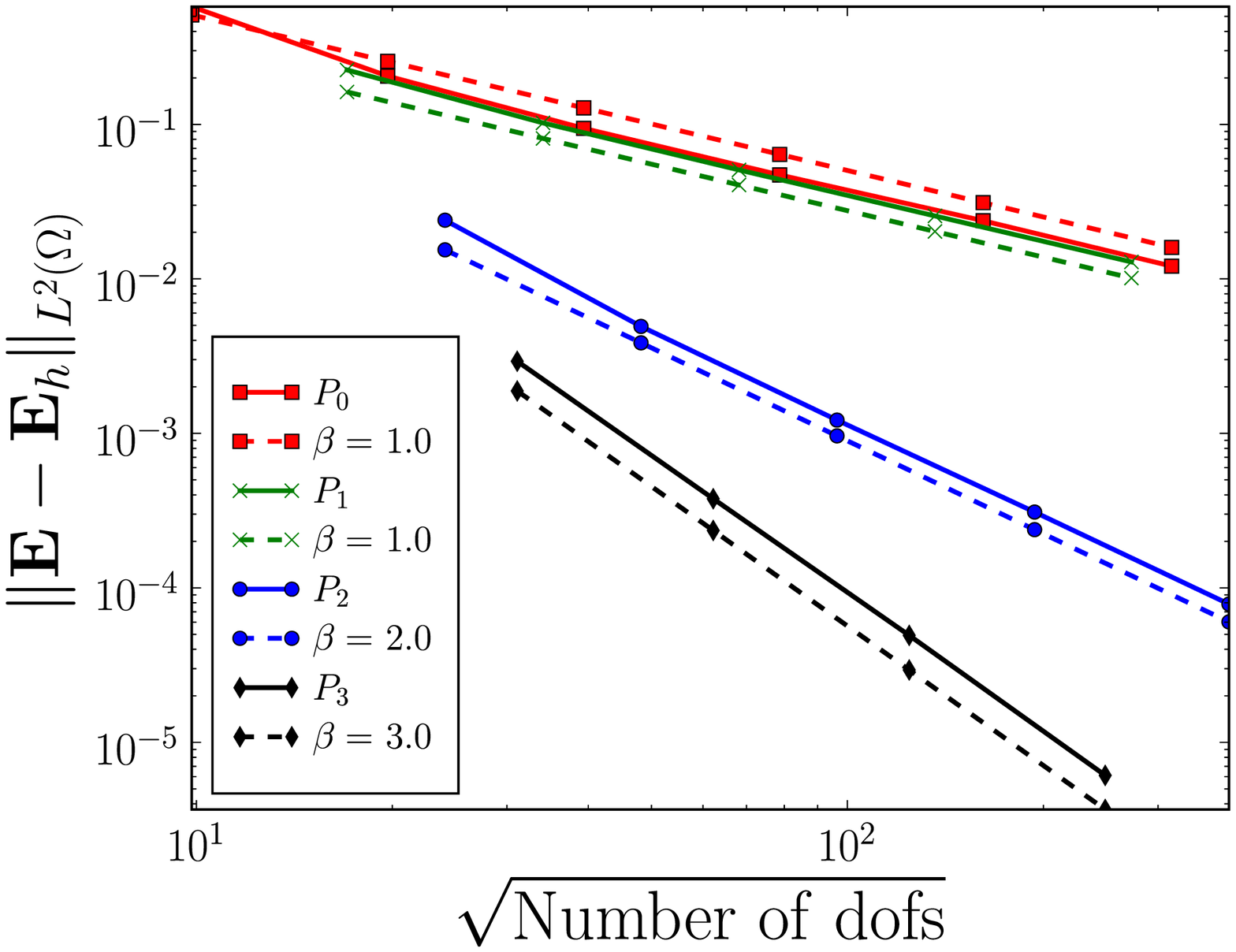} %
    \label{fig:res_c_H_ddl} } %
  \caption{Convergence results using a \emph{centered flux}.  Solid lines show
    the evolution for the whole of the numerical results and dotted lines show
    the  asymptotic  tendency, using  coefficients  $\beta$  or $\gamma$  from
    inequalities \eqref{eq:gamma_beta} estimated by a linear regression.}
  \label{fig:res_centre}
\end{figure}

\noindent {\bf Centered flux.}  Numerical convergence results in a logarithmic
scale are shown on  Figure~\ref{fig:res_centre}.  They clearly demonstrate the
interest of higher order  polynomial approximations which allow a considerable
reduction of  the number  of degrees  of freedom to  reach the  same accuracy.
Table~\ref{tab:res_centre}  summarizes  numerical  estimates (using  a  linear
regression method) of the asymptotic convergence order.

\begin{table}[!hbtp]
  \centering %
  \begin{tabular}[c]{c|c|c|c|c} %
    & $P_0$ & $P_1$ & $P_2$ & $P_3$
    \\
    \hline $\vec{E}$ & 1.0 & 1.0 & 2.0 & 3.0
    \\
    \hline $\vec{H}$ & 1.0 & 2.0 & 3.0 & 3.6
  \end{tabular}
  \caption{Numerical convergence order using a centered flux.}
  \label{tab:res_centre}
\end{table}

The method  based on  a $P_0$ approximation  (\textsl{i.e.} the  standard cell
centered finite  volume method) is  special: the convergence order  is optimal
for both fields $\vec{E}$ and $\vec{H}$,  that is, equal to $k+1$.  This could
be  the  consequence of  using  uniformly  refined  meshes, since  a  somewhat
different  behavior is obtained  for independent  meshes with  decreasing mesh
size  (see subsection  \ref{sec:maill_distincts}).  For  the  other polynomial
degrees, we  get exactly  the predicted theoretical  convergence order  in the
elliptic case for $\vec{E}$, whereas  for $\vec{H}$, this convergence order is
optimal. Therefore, in this example, the magnetic field is better approximated
than the electric field, when using the centered flux.

\begin{figure}[!hbtp] 
  \centering %
  \subfigure[$\| \vec{H} - \vec{H}_h \|_{L^2}$ against the square root of the
  number of dofs.] {
    \includegraphics[width=0.47\linewidth]{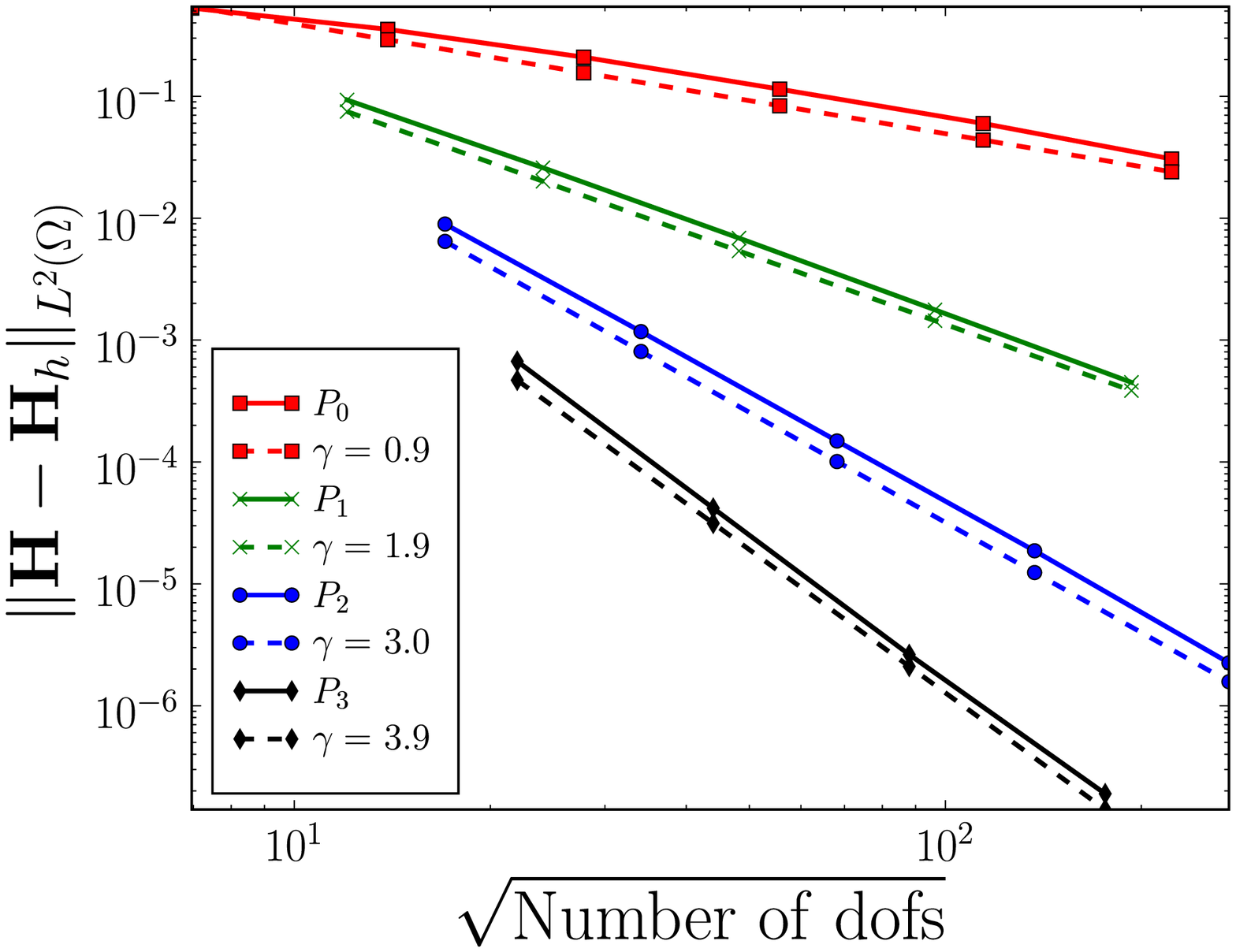} %
    \label{fig:res_d_H_pas} } %
  \hfill \subfigure[$\| \vec{E} - \vec{E}_h \|_{L^2}$ against the square root
  of the number of dofs.]  { %
    \includegraphics[width=0.47\linewidth]{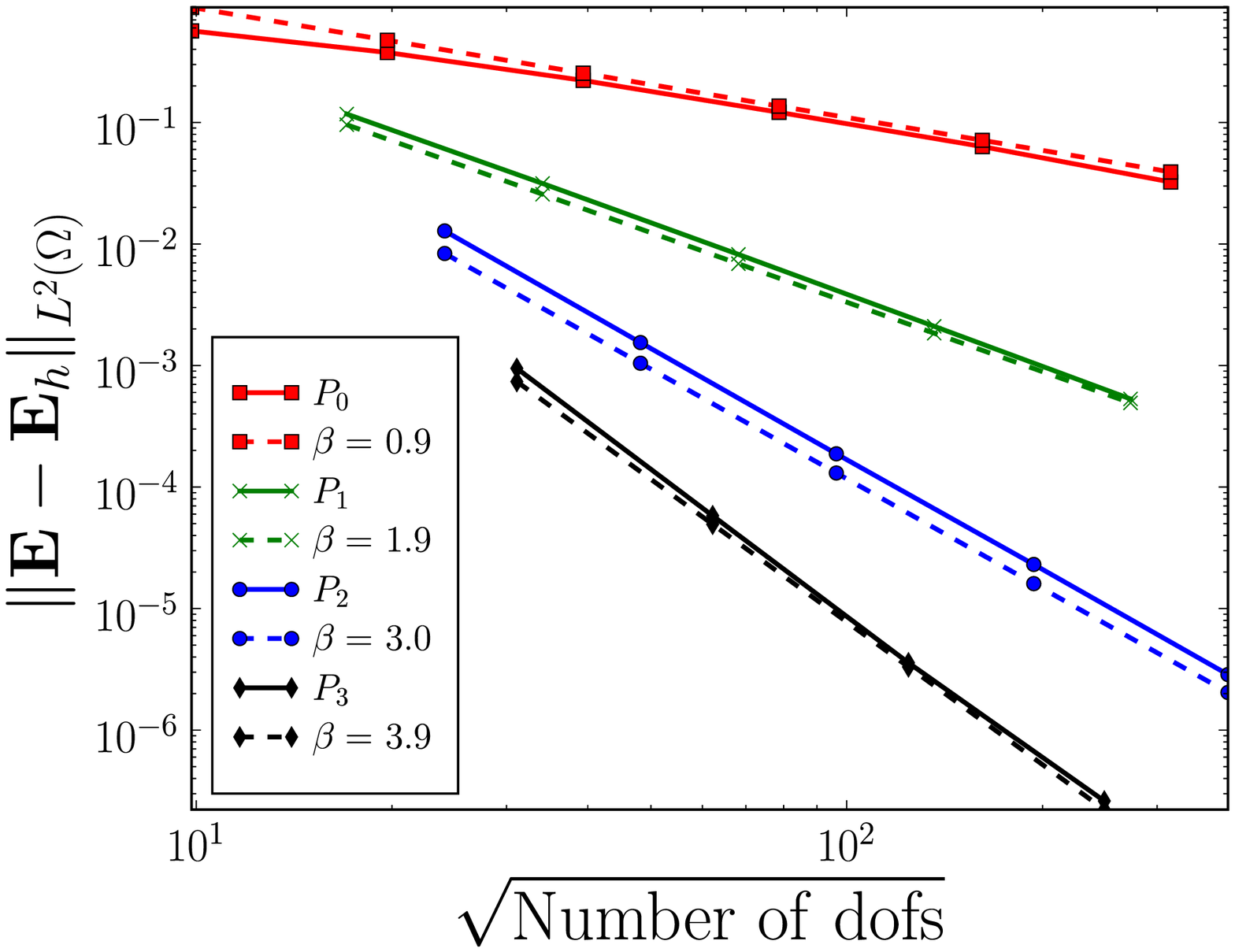} %
    \label{fig:res_d_H_ddl} } %
  \caption{Convergence results using an  \emph{upwind flux}.  Solid lines show
    the evolution for the whole of the numerical results and dotted lines show
    the  asymptotic  tendency, using  coefficients  $\beta$  or $\gamma$  from
    inequalities \eqref{eq:gamma_beta} estimated by a linear regression.}
  \label{fig:res_decentre}
\end{figure}

\noindent  {\bf  Upwind flux.}  We  used  here  the parameters  $\alpha^H_F  =
\alpha^E_F = \eta_F = 1$ for each face $F$.  Numerical convergence results are
shown on Figure~\ref{fig:res_decentre}.  Similar conclusions can be derived as
in the  centered case  except that the  convergence properties of  the methods
based on $P_0$  and $P_1$ interpolations are this  time clearly different with
respect  to  the  centered  case.   The  asymptotic  convergence  orders  (see
Table~\ref{tab:res_decentre}) are  similar for  both fields and  correspond to
the  theory for  the elliptic  Maxwell equations.  The convergence  is optimal
except for the case $P_0$, but nevertheless we are still above the theoretical
estimates.

\begin{table}[!hbtp] %
  \centering %
  \begin{tabular}[c]{c|c|c|c|c} %
    & $P_0$ & $P_1$ & $P_2$ & $P_3$
    \\
    \hline $\vec{E}$ & 0.9 & 1.9 & 3.0 & 3.9
    \\
    \hline $\vec{H}$ & 0.9 & 1.9 & 3.0 & 3.9
  \end{tabular}
  \caption{Numerical convergence order using an upwind flux.}
  \label{tab:res_decentre}
\end{table}

\begin{figure}[!hbtp]
  \centering %
  \subfigure[$\| \vec{H} - \vec{H}_h \|_{L^2}$ against the square root of the
  number of dofs.] {
    \includegraphics[width=0.47\linewidth]{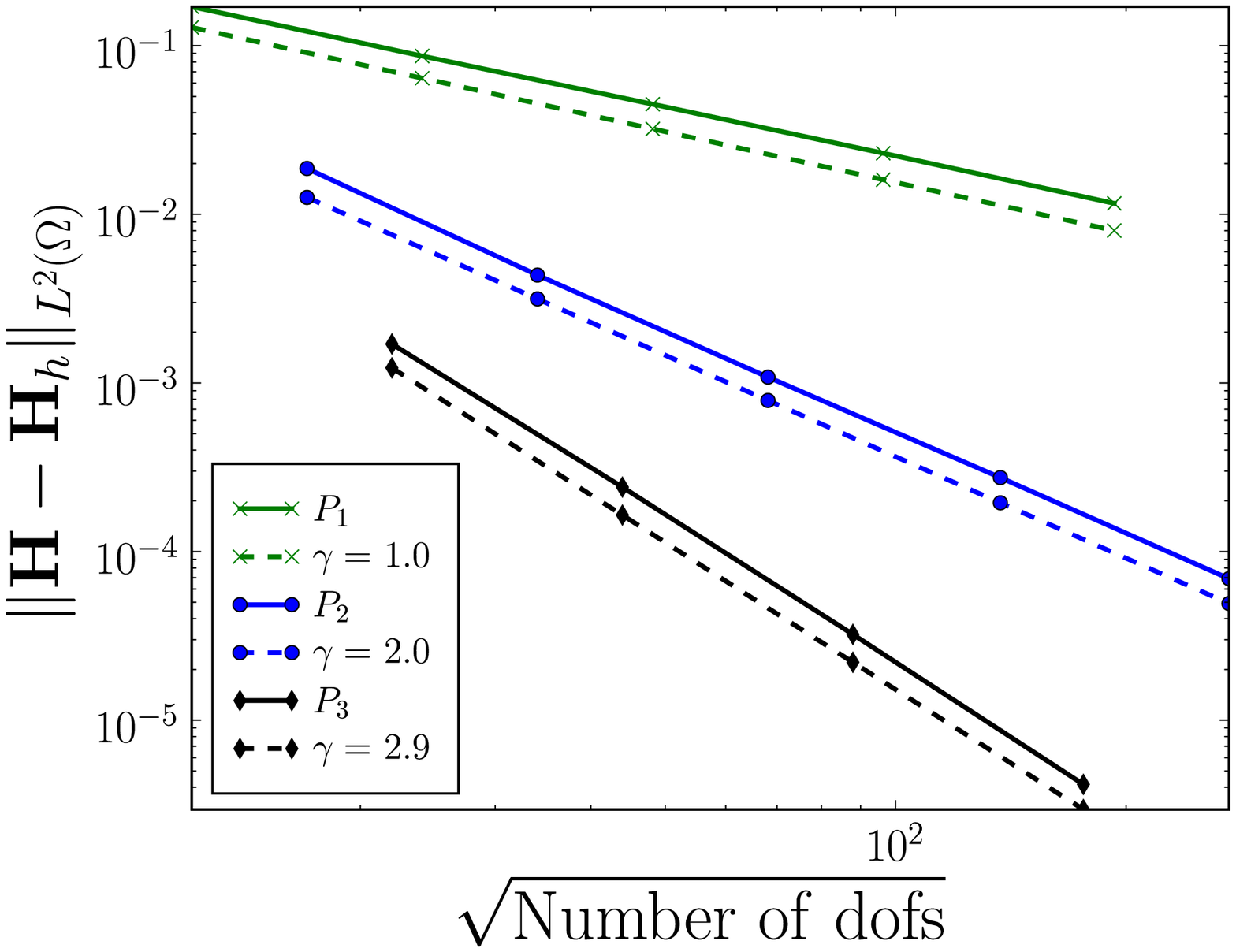} %
    \label{fig:res_p_H_pas} } %
  \hfill \subfigure[$\| \vec{E} - \vec{E}_h \|_{L^2}$ against the square root
  of the number of dofs.]  { %
    \includegraphics[width=0.47\linewidth]{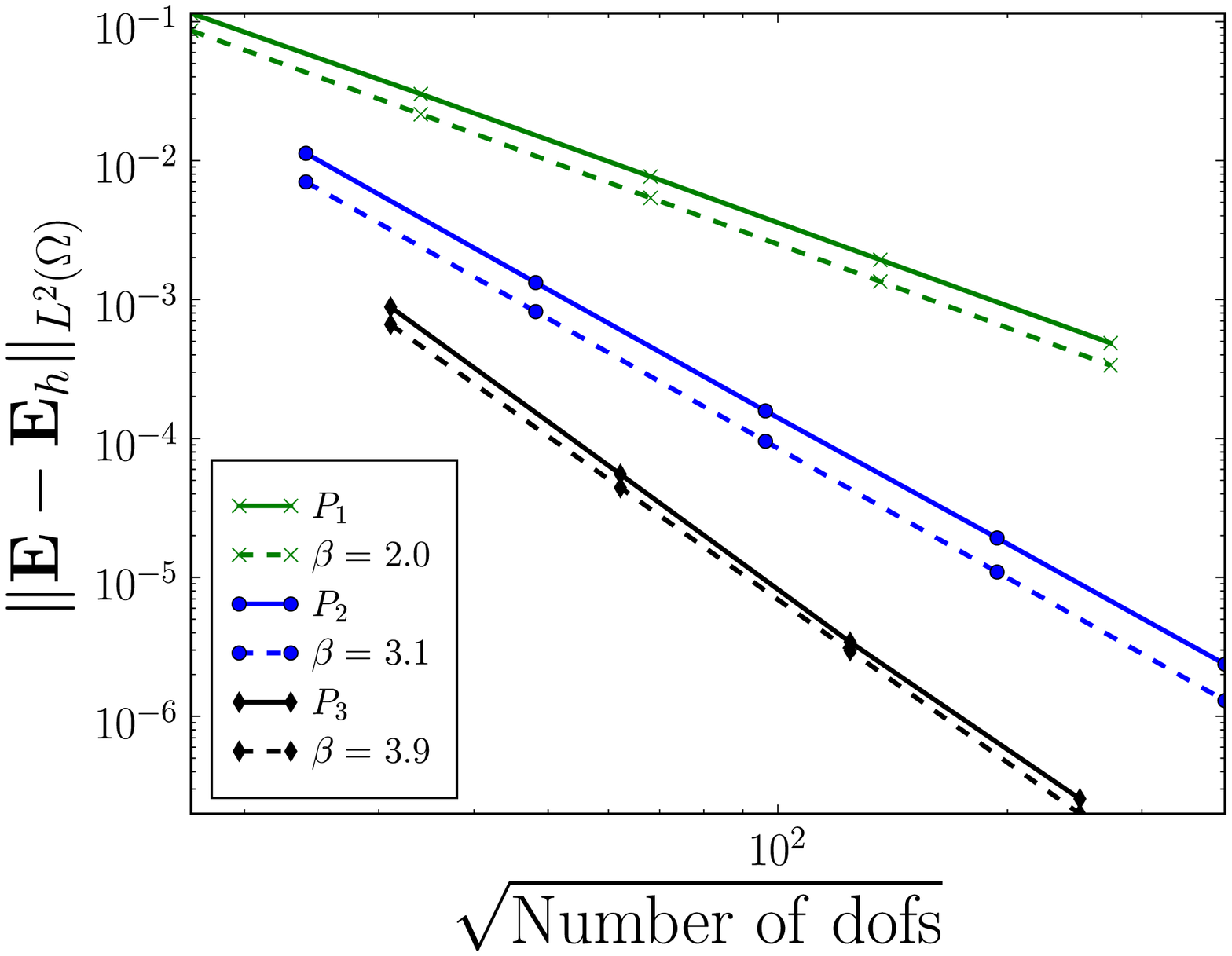} %
    \label{fig:res_p_H_ddl} } %
  \caption{Convergence  results using  a \emph{penalized  flux  on $\vec{E}$}.
    Solid lines show the evolution for  the whole of the numerical results and
    dotted lines  show the asymptotic tendency, using  coefficients $\beta$ or
    $\gamma$  from inequalities  \eqref{eq:gamma_beta} estimated  by  a linear
    regression.}
  \label{fig:res_penal}
\end{figure}

\noindent {\bf Penalized flux on $\vec{E}$.}  We set $\tau_F = \eta_F = 1$ for
each  face  $F$.  Results   are  shown  on  Figure~\ref{fig:res_penal}.  Table
\ref{tab:res_penal}  summarizes  the  numerical  estimates of  the  asymptotic
convergence order. Besides the expected lack of convergence in the case $P_0$,
we can notice for all the other cases ($(P_k)_{k>0}$) a complementary behavior
with  respect  to  the centered  flux,  since  this  time  we get  an  optimal
convergence rate for $\vec{E}$, but not for $\vec{H}$.

\begin{table}[!hbtp]
  \centering %
  \begin{tabular}[c]{c|c|c|c|c} %
    & $P_0$ & $P_1$ & $P_2$  & $P_3$
    \\
    \hline $\vec{E}$ & X & 2.0 & 3.1 & 3.9
    \\
    \hline $\vec{H}$ & X & 1.0 & 2.0 & 2.9
  \end{tabular}
  \caption{Numerical convergence order using a penalized flux on $\vec{E}$.}
  \label{tab:res_penal}
\end{table}

\subsection{Convergence behavior using independent meshes}
\label{sec:maill_distincts}

\begin{figure}[!hbtp] %
  \centering %
  \subfigure[$\| \vec{H} - \vec{H}_h \|_{L^2}$ against the mesh size $h$.  ]{
    \label{fig:CompHtu}
    \includegraphics[width=0.47\linewidth]{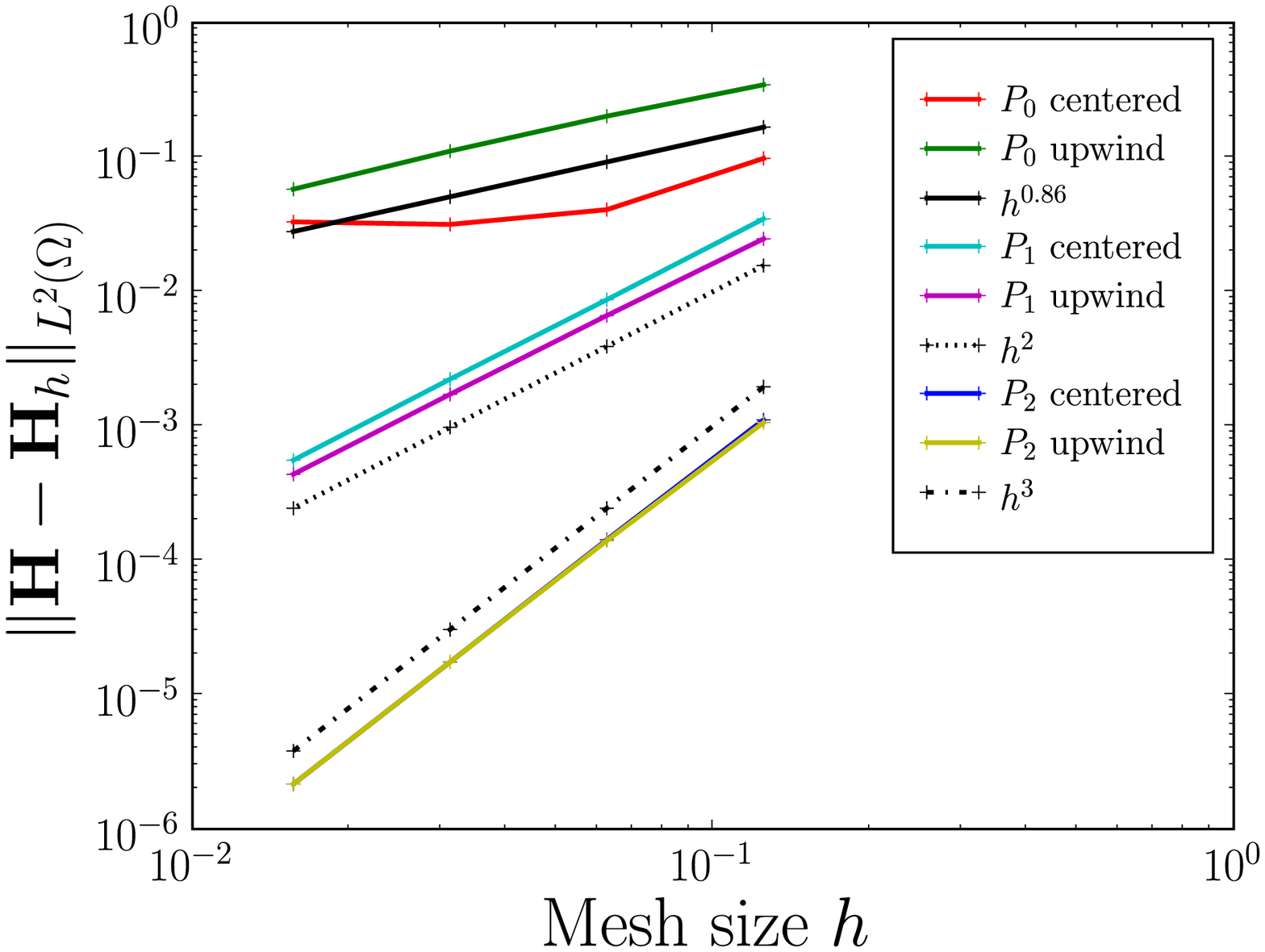} } %
  \hfill \subfigure[$\| \vec{E} - \vec{E}_h \|_{L^2}$ against the mesh size
  $h$.  ]{ %
    \label{fig:CompEtu}
    \includegraphics[width=0.47\linewidth]{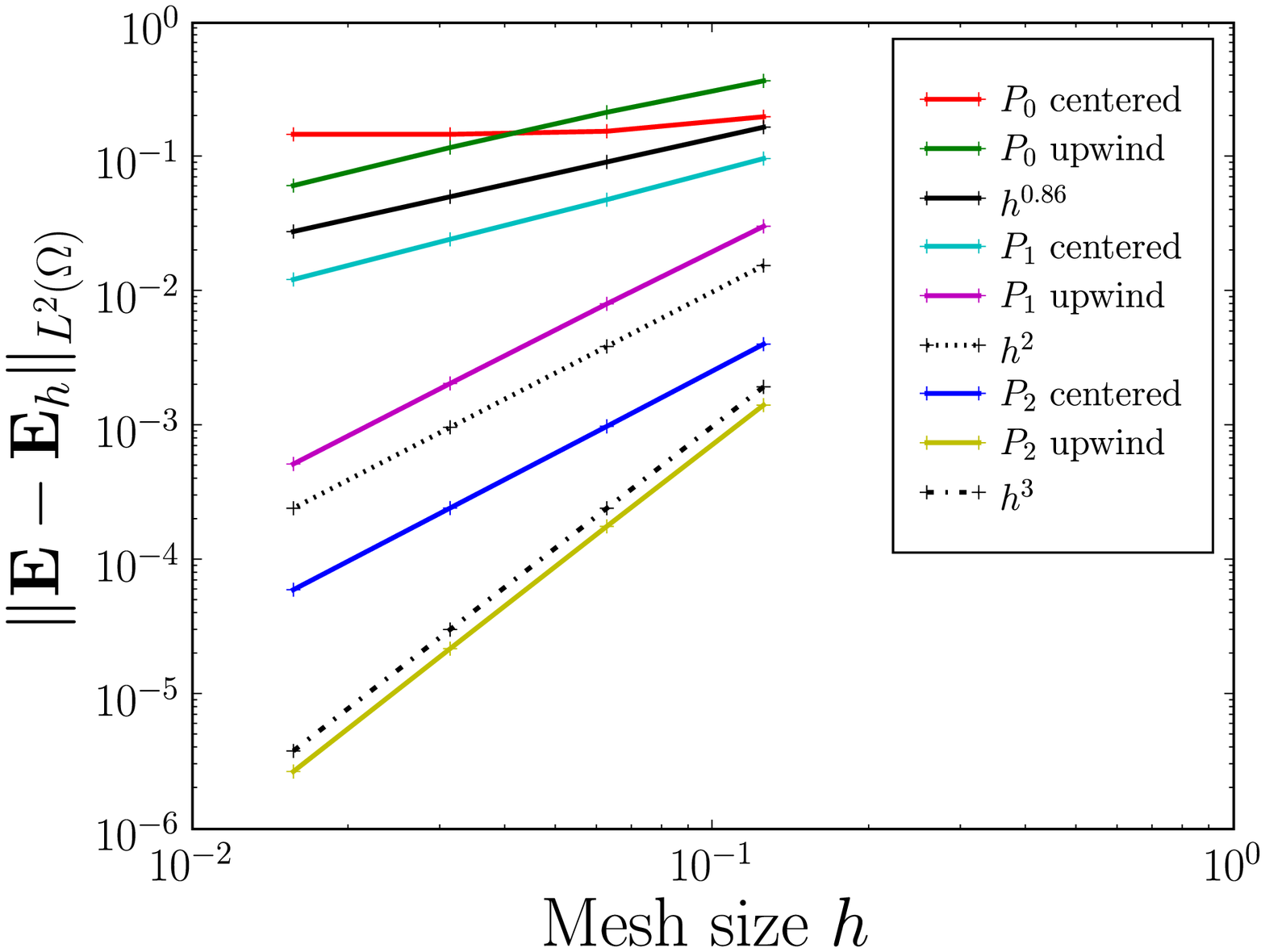} }
  \caption{Comparison of the convergence results between \emph{centered flux}
    and \emph{upwind flux}.}
  \label{fig:compare_tu}
\end{figure}

On Figure~\ref{fig:compare_tu}, we compare  the evolution of the $L^2$-norm of
the    error   with    the   mesh    size    $h$   by    using   the    meshes
$(\Trond_i)_{i=1,\ldots,4}$,  for both  a centered  flux and  an  upwind flux,
Figure~\ref{fig:CompEtu} corresponds to the error for the field $\vec E$ while
Figure~\ref{fig:CompHtu} corresponds to the error for the field $\vec H$.  The
results for the upwind flux are  the same as for the uniformly refined meshes.
For the  centered flux, note the lack  of convergence for the  case $P_0$. For
all the other  cases the results remain the same as  for the uniformly refined
meshes.

It is already  known for time-domain problems that  the centered flux combined
to  a   leap-frog  time  integration  scheme  results   in  a  non-dissipative
discontinuous Galerkin method (a mandatory feature for long time computations,
see  \cite{fezoui05}).  As far  as time-harmonic  problems are  concerned, the
previous results show that the  upwind flux has better convergence properties.
Nevertheless, the  centered flux remains  less expensive both  for time-domain
and  time-harmonic problems (arithmetic  operations and  memory requirements).

\subsection{Numerical comparisons on a less trivial problem}
\label{sec:essainum_2}

\begin{figure}[!hbtp]
  \centering %
  \subfigure[First mesh. $h_{\max} = 0.32$]{%
    \includegraphics[width=0.3\linewidth]{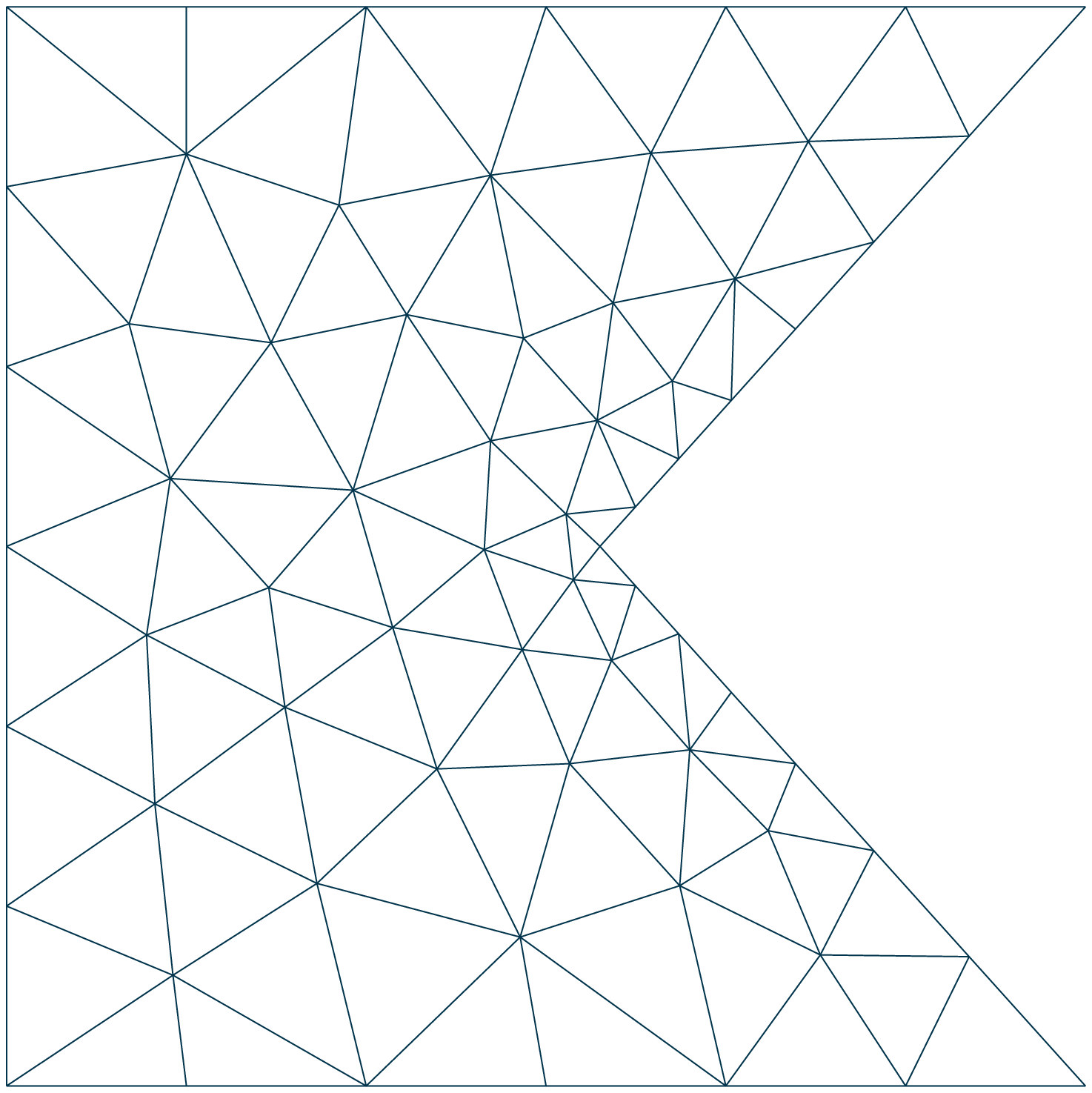} %
  } \hfill%
  \subfigure[Second mesh. $h_{\max} = 0.16$]{%
    \includegraphics[width=0.3\linewidth]{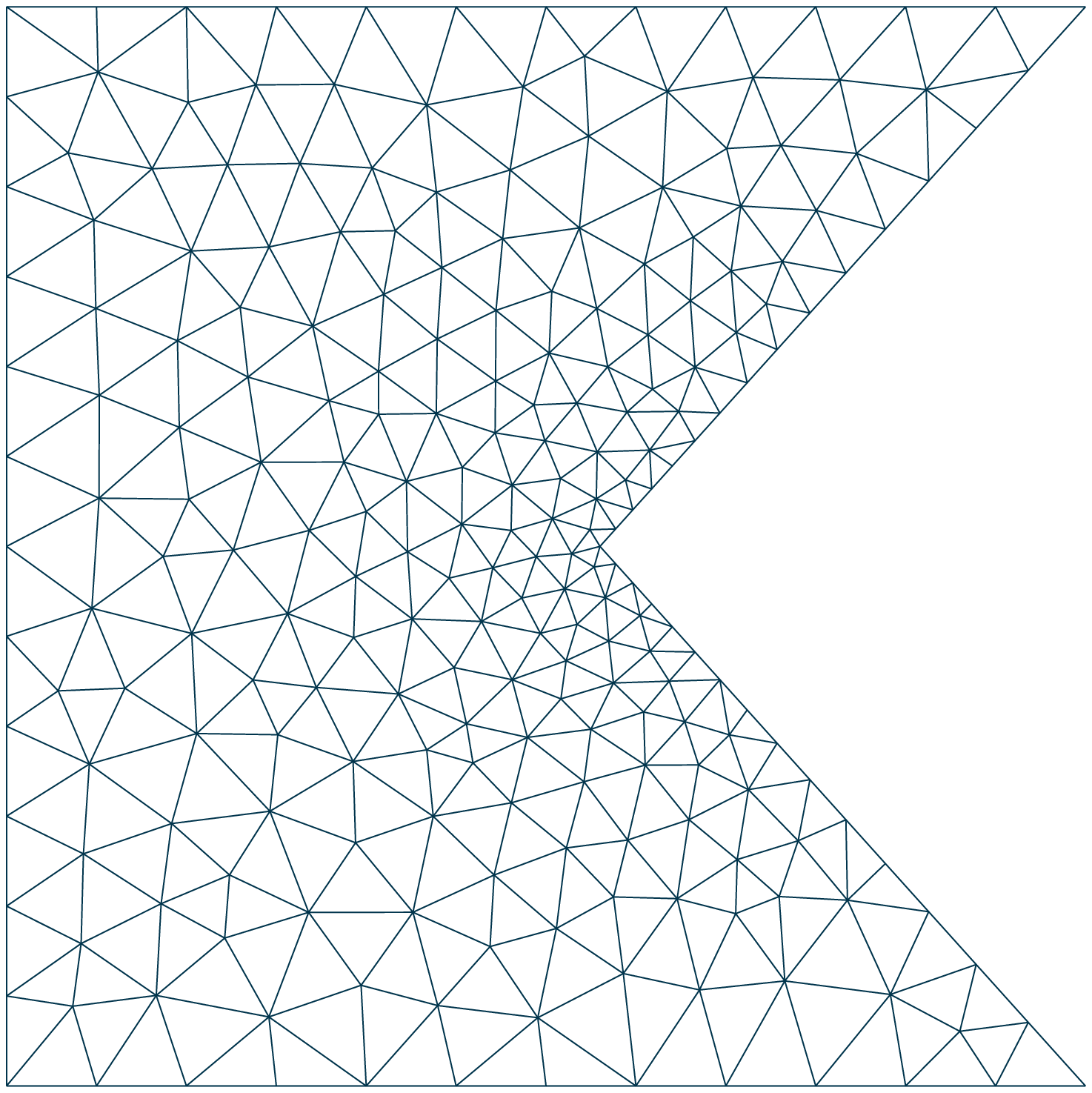} %
  } \hfill%
  \subfigure[Third mesh. $h_{\max} = 0.32$]{%
    \includegraphics[width=0.3\linewidth]{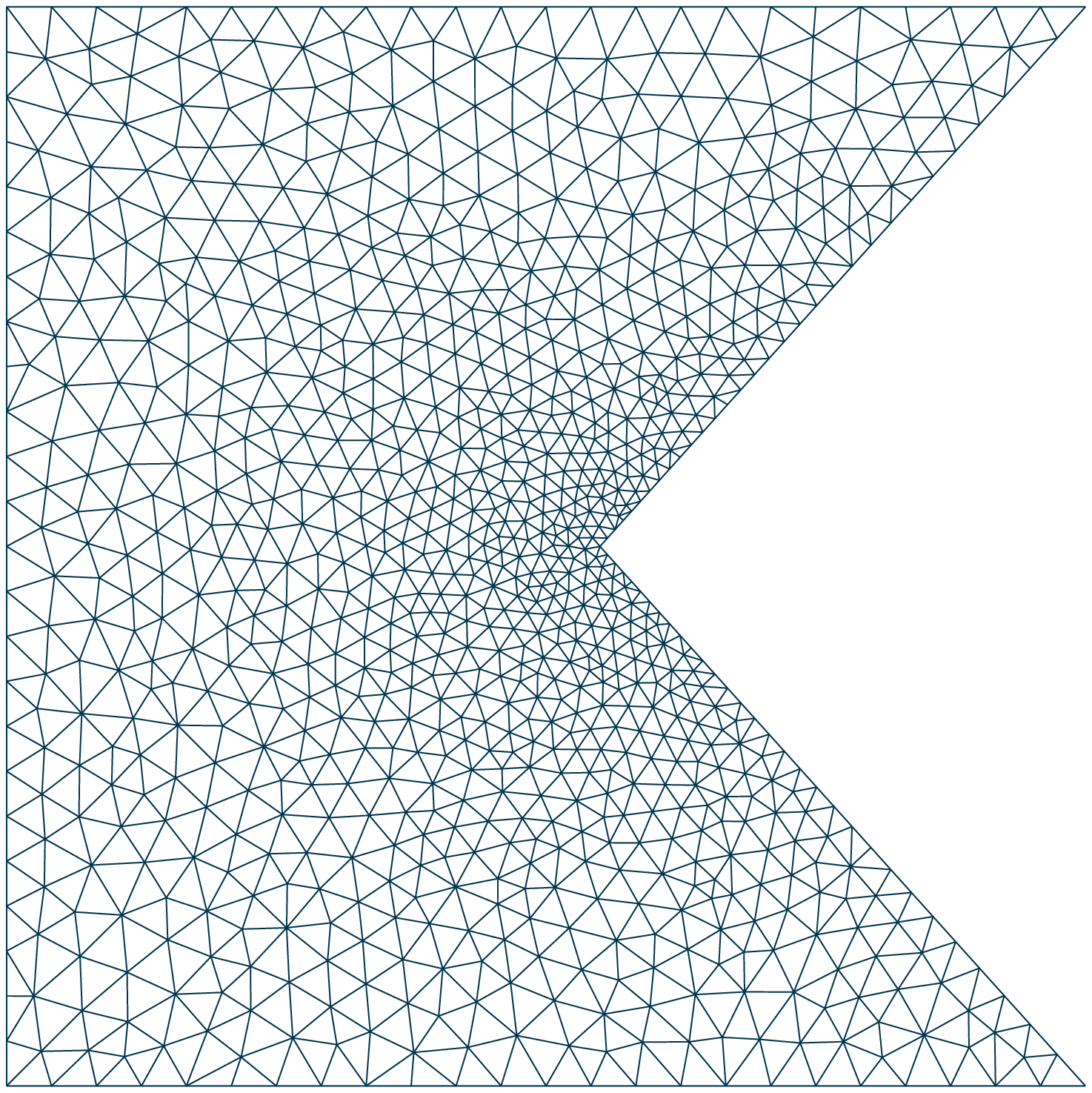} %
  }
  \caption{Three first meshes used for the second example.}
  \label{fig:num2_meshes}
\end{figure}

The  domain is  the square  $[-1; 1]^2$  where we  have suppressed  a  part by
inserting   a  point   of  coordinates   $(0.1,  0)$   at  it   is   shown  on
Figure~\ref{fig:num2_meshes}.   The properties  $\epsilon_r$  and $\mu_r$  are
still  homogeneous  and  equal  to one.Appropriate  non-homogeneous  Dirichlet
boundary conditions  are enforced on  the boundary of  the domain in  order to
obtain $\vec E = (\sin (2\pi y), \sin (2\pi x))^t$ as the solution.

The    mesh   is    not    fully    homogeneous   as    it    is   shown    on
Figure~\ref{fig:num2_meshes};  it is  slightly  denser next  to  the point  of
coordinates   $(0.1,  0)$.    Independent  meshes   have  been   used   as  in
Subsection~\eqref{sec:maill_distincts}.

\begin{figure}[!hbtp]
  \centering
  \subfigure[Evolution of the $L^2$-norm of the error for the $\vec E$
  field.]{%
    \includegraphics[width=0.45\linewidth]{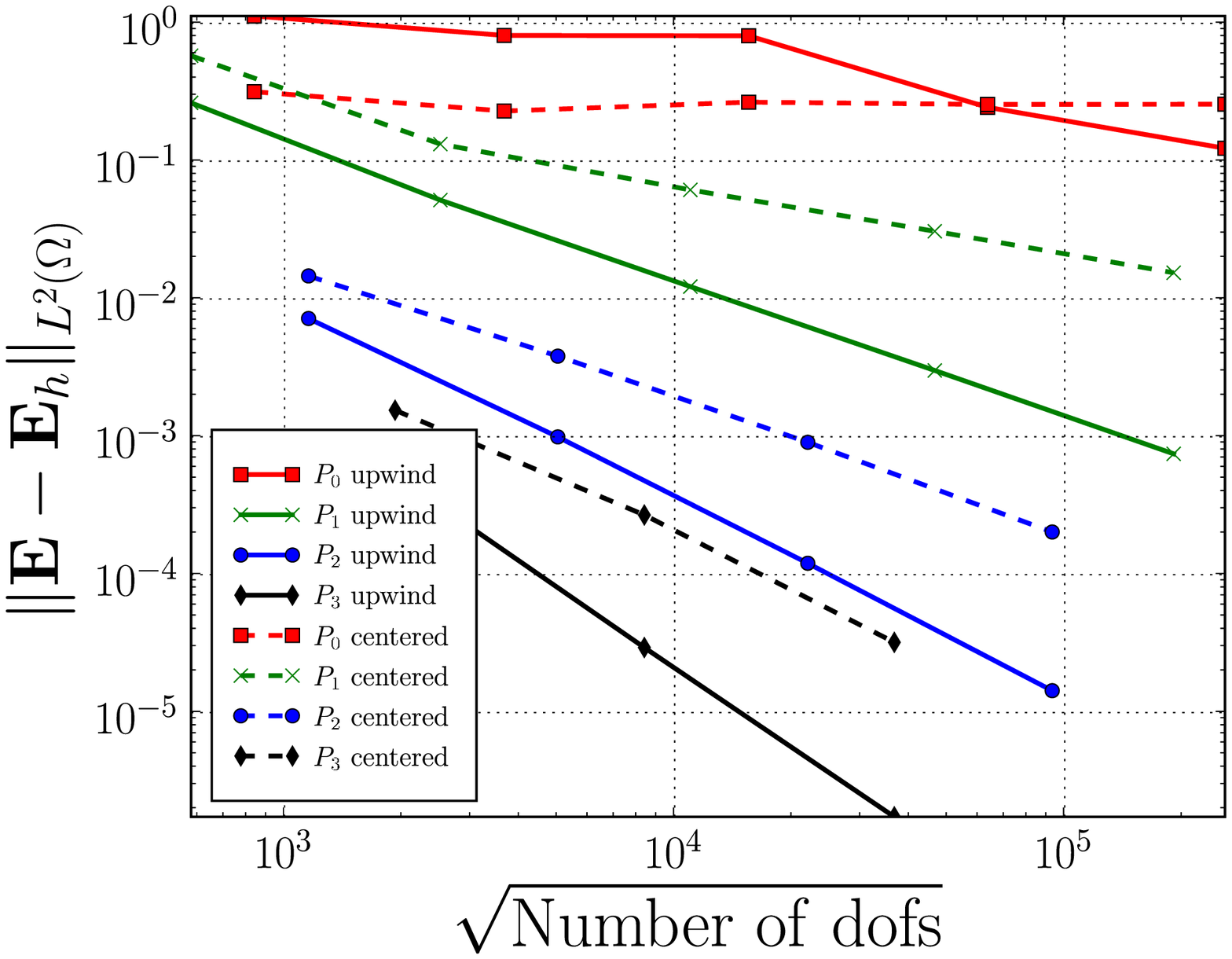}%
  } \hfill%
  \subfigure[Evolution of the $L^2$-norm of the error for the $\vec H$
  field.]{%
    \includegraphics[width=0.45\linewidth]{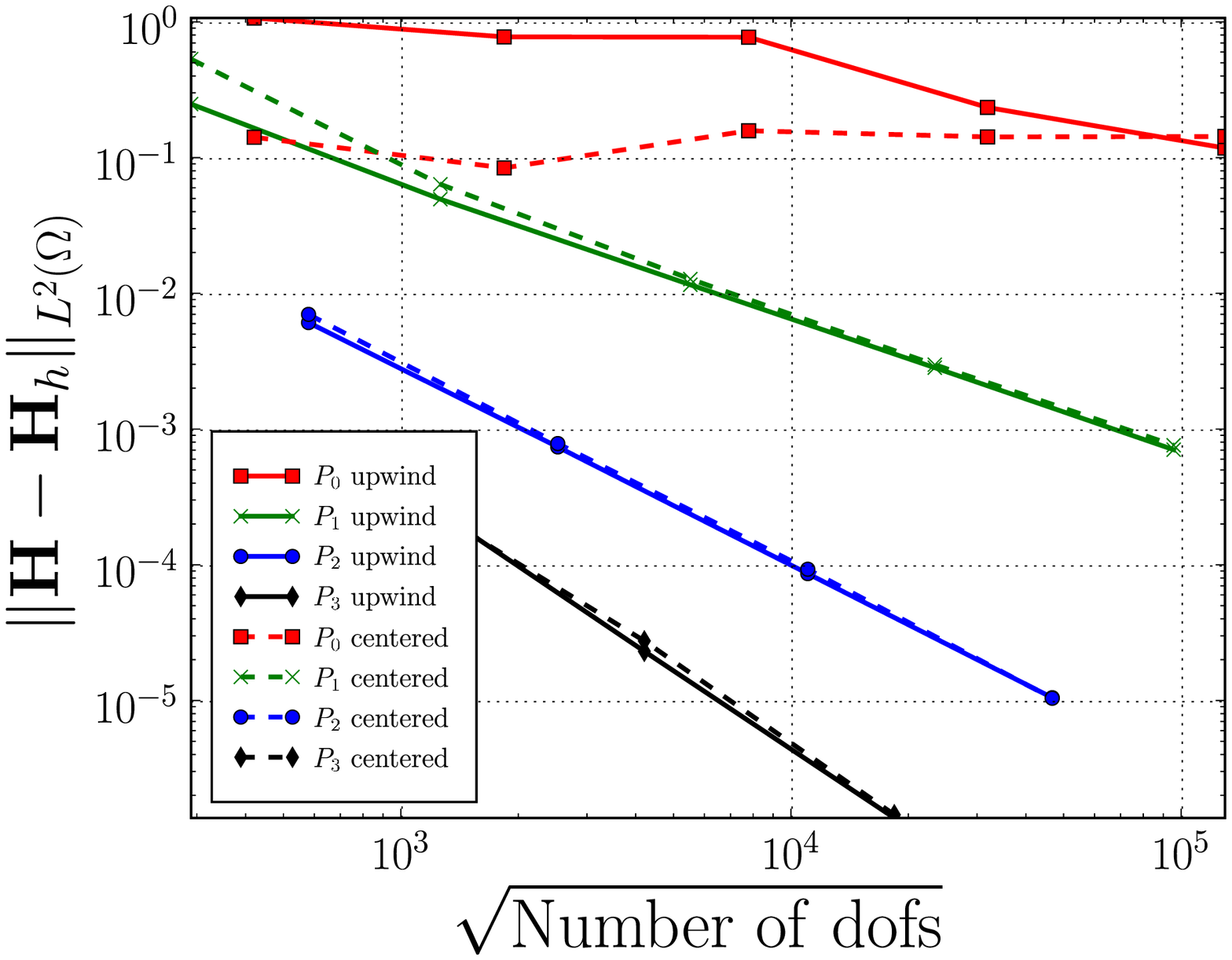}%
  }%
  \caption{Evolution of the $L^2$-norm of the error against the square root of
    the number of degrees of freedom (dofs).}
  \label{fig:num2_compare}
\end{figure}

The  results shown on  Figure~\ref{fig:num2_compare} are  in a  full agreement
with those  obtained in the case of  the plane wave and  independent meshes at
Subsection \ref{sec:maill_distincts}.

\end{document}